\documentclass[letterpaper]{amsart}

\newtheorem{theorem}{Theorem}[section]
\newtheorem{lemma}[theorem]{Lemma}

\newtheorem{proposition}[theorem]{Proposition}
\newtheorem{corollary}[theorem]{Corollary}

\theoremstyle{definition}
\newtheorem{definition}[theorem]{Definition}

\theoremstyle{remark}
\newtheorem{remark}[theorem]{Remark}

\usepackage{graphicx}
\usepackage{color}
\usepackage{amsmath}
\usepackage{amsfonts}
\usepackage{amssymb}
\usepackage{epsfig}
\usepackage{rotating}
\usepackage{graphics}
\newcommand{\be}{\begin{equation}}

\newcommand{\ee}{\end{equation}}

\newcommand{\ts}{\textstyle}

\newcommand{\la}{\lambda}

\newcommand{\dz}{\wedge}

\newcommand{\ba}{\begin{array}}

\newcommand{\ea}{\end{array}}

\newcommand{\beq}{\begin{eqnarray}}

\newcommand{\eeq}{\end{eqnarray}}

\newtheorem{lm}{lemma}

\newtheorem{thee}{theorem}

\newtheorem{proo}{proposition}

\newtheorem{co}{corollary}

\newtheorem{rem}{remark}

\newtheorem{deff}{definition}

\newcommand{\bd}{\begin{deff}}

\newcommand{\ed}{\end{deff}}

\newcommand{\bl}{\begin{lm}}

\newcommand{\el}{\end{lm}}

\newcommand{\bp}{\begin{proo}}

\newcommand{\ep}{\end{proo}}

\newcommand{\bt}{\begin{thee}}

\newcommand{\et}{\end{thee}}

\newcommand{\bc}{\begin{co}}

\newcommand{\ec}{\end{co}}

\newcommand{\brm}{\begin{rem}}

\newcommand{\erm}{\end{rem}}

\newcommand{\der}{{\rm d}}

\usepackage{t1enc}
\def\frak{\mathfrak}

\newcommand{\newc}{\newcommand}

\let\ccdot\cdot
\def\cdot{\hbox to 2.5pt{\hss$\ccdot$\hss}}

\newc{\aR}{\mbox{\boldmath{$ R$}}}
\newc{\aS}{\mbox{\boldmath{$ S$}}}
\newc{\aT}{\mbox{\boldmath{$ T$}}}
\newc{\aW}{\mbox{\boldmath{$ W$}}}

\newc{\aK}{\mbox{\boldmath{$ K$}}}
\newc{\aL}{\mbox{\boldmath{$ L$}}}

\newcommand{\bbP}{\mathbb{P}}
\newcommand{\bbC}{\mathbb{C}}

\usepackage{amssymb}
\usepackage{amscd}

\newcommand{\hook}{\raisebox{-0.35ex}{\makebox[0.6em][r]
{\scriptsize $-$}}\hspace{-0.15em}\raisebox{0.25ex}{\makebox[0.4em][l]{\tiny
 $|$}}}

\newcommand{\bma}{\begin{pmatrix}}
\newcommand{\ema}{\end{pmatrix}}

\newc{\obstrn}[2]{B^{#1}_{#2}}

\newcommand{\rpl}                         
{\mbox{$
\begin{picture}(12.7,8)(-.5,-1)
\put(0,0.2){$+$}
\put(4.2,2.8){\oval(8,8)[r]}
\end{picture}$}}

\newcommand{\lpl}                         
{\mbox{$
\begin{picture}(12.7,8)(-.5,-1)
\put(2,0.2){$+$}
\put(6.2,2.8){\oval(8,8)[l]}
\end{picture}$}}

\usepackage{ifthen}

\newc{\tensor}[1]{#1}
\newc{\Mvariable}[1]{\mbox{#1}}
\newc{\down}[1]{{}_{#1}}
\newc{\up}[1]{{}^{#1}}

\newc{\JulyStrut}{\rule{0mm}{6mm}}
\newc{\midtenPan}{\mbox{\sf S}}
\newc{\midten}{\mbox{\sf T}}
\newc{\midtenEi}{\mbox{\sf U}}
\newc{\ATen}{\mbox{\sf E}}
\newc{\BTen}{\mbox{\sf F}}
\newc{\CTen}{\mbox{\sf G}}

\def\sideremark#1{\ifvmode\leavevmode\fi\vadjust{\vbox to0pt{\vss
 \hbox to 0pt{\hskip\hsize\hskip1em
 \vbox{\hsize3cm\tiny\raggedright\pretolerance10000
 \noindent #1\hfill}\hss}\vbox to8pt{\vfil}\vss}}}%

\numberwithin{equation}{section}

\newcounter{romenumi}
\newcommand{\labelromenumi}{(\roman{romenumi})}

\begin{document}
\title{Einstein's equations and the embedding of 3-dimensional CR manifolds}
\vskip 1.truecm
\author{C. Denson Hill} \address{Department of Mathematics, Stony
  Brook University, Stony Brook, N.Y. 11794, USA}
\email{dhill@math.sunysb.edu}  \thanks{This research was supported in part by
the Polish grant 1 P03B 07529}
\author{Jerzy Lewandowski} \address{Instytut Fizyki Teoretycznej,
Uniwersytet Warszawski, ul. Hoza 69, Warszawa, Poland}
\email{lewand@fuw.edu.pl}
\author{Pawe\l~ Nurowski} \address{Instytut Fizyki Teoretycznej,
Uniwersytet Warszawski, ul. Hoza 69, Warszawa, Poland}
\email{nurowski@fuw.edu.pl}

\date{\today}

\begin{abstract}
We prove several theorems concerning the connection between the local
CR embeddability of 3-dimensional CR manifolds, and the existence of
algebraically special Maxwell and gravitational fields. We reduce the Einstein equations for spacetimes associated with such
fields to a system of 
CR invariant equations on a 3-dimensional CR manifold defined by the fields. 
Using the reduced Einstein 
equations we construct two independent CR functions for the
corresponding CR manifold. We also point out that the Einstein
equations, imposed on spacetimes associated with a 3-dimensional 
CR manifold, imply that the spacetime metric, after an appropriate
rescaling, becomes 
well defined on a circle bundle over the CR manifold. The circle
bundle itself emerges as a consequence of Einstein's equations.
\vskip5pt\centerline{\small\textbf{MSC classification}: 32V10, 32V30,
  83C15, 83C20}\vskip15pt
\end{abstract}
\maketitle
\tableofcontents
\newcommand{\bbS}{\mathbb{S}}
\newcommand{\bbR}{\mathbb{R}}
\newcommand{\sog}{\mathbf{SO}}
\newcommand{\slg}{\mathbf{SL}}
\newcommand{\og}{\mathbf{O}}
\newcommand{\soa}{\frak{so}}
\newcommand{\sla}{\frak{sl}}
\newcommand{\sua}{\frak{su}}
\newcommand{\dr}{\mathrm{d}}
\newcommand{\sug}{\mathbf{SU}}
\newcommand{\gat}{\tilde{\gamma}}
\newcommand{\Gat}{\tilde{\Gamma}}
\newcommand{\thet}{\tilde{\theta}}
\newcommand{\Thet}{\tilde{T}}
\newcommand{\rt}{\tilde{r}}
\newcommand{\st}{\sqrt{3}}
\newcommand{\kat}{\tilde{\kappa}}
\newcommand{\kz}{{K^{{~}^{\hskip-3.1mm\circ}}}}
\newcommand{\bv}{{\bf v}}
\newcommand{\di}{{\rm div}}
\newcommand{\curl}{{\rm curl}}
\newcommand{\cs}{(M,{\rm T}^{1,0})}
\section{Introduction}
\subsection{3-dimensional CR structures}
Let $M$ be an open set in $\bbR^3$, with \emph{real}
coordinates $(x,y,z)$. We consider a complex vector field 
\be
Z=a\frac{\partial}{\partial x}+b\frac{\partial}{\partial y}+c\frac{\partial}{\partial z}.\label{dbar}
\ee
Here the coefficients 
$$a=a(x,y,z),\quad\quad b=b(x,y,z),\quad\quad
c=c(x,y,z)$$  
are \emph{complex valued} functions on $M$. 
We assume that that the vector field $Z$ and its complex
conjugate $\bar{Z}$ are \emph{linearly independent} at each point of
$M$. The vector field $Z$ spans a 1-dimensional complex
distribution ${\rm T}^{1,0}$ in $\bbC{\rm T}M$. By definition the pair
$\cs$ is 
an abstract \emph{3-dimensional CR manifold}. The CR structure $\cs$
on $M$  
is said to be of class $C^k$ iff the coefficients 
$a=a(x,y,z)$, $b=b(x,y,z)$ and $c=c(x,y,z)$ are of differentiability 
class $C^k$. 

If the vector fields $Z$,
$\bar{Z}$ and $[Z,\bar{Z}]$
are linearly independent at each point of $M$ then the CR
structure $\cs$ is called \emph{strictly pseudoconvex}. When
$[Z,\bar{Z}]$ is linearly dependent on $Z$ and $\bar{Z}$ at a point
$p$, the \emph{Levi form} of the structure $\cs$ vanishes 
at $p$.

Natural examples of 3-dimensional CR manifolds are hypersurfaces
$M$ of one real codimension embedded in $\bbC^2$. They 
acquire a CR structure ${\rm T}^{1,0}M$ from the ambient space by ${\rm
  T}^{1,0}M=\{X-iJX~|~ X\in {\rm
  T}M\cap J({\rm T}M)\}$, where $J$ is the
standard complex structure on $\bbC^2$. 

Given an abstract CR manifold $\cs$ one asks if there exists a
local diffeomorphism  
$$\iota: M\to \iota(M)\subset\bbC^2$$
such that
$$\iota_* {\rm T}^{1,0}= {\rm
  T}^{1,0}\iota(M).$$
Such a diffeomorphism is called a \emph{local CR embedding} of $\cs$ 
into $\bbC^2$.  

The question of whether or not  a given 3-dimensional CR manifold $\cs$ 
can be locally CR embedded as a hypersurface in $\bbC^2$ is related to
the problem of the local existence of solutions to the \emph{linear}
partial differential equation 
\be
\bar{Z}\zeta = 0,\label{hl}
\ee
where
$$\bar{Z}=\bar{a}\frac{\partial}{\partial
  x}+\bar{b}\frac{\partial}{\partial
  y}+\bar{c}\frac{\partial}{\partial z}$$
is a section of the bundle $\bar{\rm T}^{1,0}={\rm T}^{0,1}$. The
  solution $\zeta=\zeta(x,y,z)$ is
called a \emph{CR function}. Of course any holomorphic function
$h=h(\zeta)$ of a CR function is a CR function. When speaking
about different CR functions, we will mean only those which are not
functionally related to each other by a holomorphic function in this sense. 

If the CR structure $\cs$ on $M$ is \emph{real analytic}, then
by the Cauchy-Kowalewska theorem, equation (\ref{hl}) locally admits
\emph{two} functionally independent CR functions $$\zeta=\zeta(x,y,z),\quad\quad
\eta=\eta(x,y,z)$$ 
such that
$$\der\zeta\dz\der\eta\neq 0.$$ Then the CR manifold is locally CR
embeddable with the local CR embedding being real analytic, and given
by \cite{ah,Le1} 
$$\iota: M\ni (x,y,z)\mapsto
(\zeta,\eta)=(\zeta(x,y,z),\eta(x,y,z))\in \bbC^2.$$

The situation in which (\ref{hl}) has two functionally independent solutions may also happen when the CR structure on $M$ is sufficiently smooth but \emph{not} real analytic. 
The important thing is that the requirement that (\ref{hl}) locally 
admits two functionally independent CR functions is \emph{equivalent} to the local embeddability of $M$. The real analyticity is not needed for this equivalence to hold.   

It turns out that if one abandons the assumption about the real
analyticity of $\cs$, e.g. if one assumes that
$\cs$ is only of class $C^\infty$, then the equation
(\ref{hl}) \emph{may have no other local solutions than the trivial ones}
$\zeta=const$. This remarkable result is due to Louis Nirenberg
\cite{ni1,ni2}, where he gave the first example of a $Z$ having
$C^\infty$ coefficients such that (\ref{hl}) has only constant local
solutions. As a consequence such abstract CR manifolds are \emph{not} locally
CR embeddable as hypersurfaces in $\bbC^2$.

There are a few situations in which a $C^\infty$ 3-dimensional
abstract CR manifold \emph{is} locally embeddable. Among them are the
\emph{Levi flat} structures (those whose Levi form vanishes
identically in a neighborhood), as well as those which admit a local
\emph{symmetry} (i.e. a local real
vector field $X$ such that
$[X,Z]\dz Z=0$, in a neighborhood). It is
also possible to find $C^\infty$ 3-dimensional abstract CR structures
on $M$ which have one CR function $\zeta$, with $\der\zeta\neq 0$, but such
that any other local CR function is functionally dependent on $\zeta$;
hence such structures are also not locally CR embeddable. 

However the surprising and remarkable thing about $C^\infty$
3-dimensional CR structures is that \emph{the generic situation is like the
example of Nirenberg}: The complex vector fields $Z$, having only
locally constant solutions, \emph{are dense} (in the $C^\infty$
topology). This is a consequence of a Baire category argument. 

Given a strictly pseudoconvex structure $\cs$ on $M$ which is locally CR
embeddable, an arbitrary small randomly chosen $C^\infty$ perturbation
of it will no longer be locally CR embeddable, see \cite{ja}. 

Consequently, for smooth 3-dimensional CR structures, those which are
locally embeddable form a \emph{thin set} (in the sense of Baire
category) in the space of all such structures. So we have the question: What sort of conditions (aside
from real analyticity) can be imposed on a sufficiently smooth
3-dimensional CR structure in order to guarantee local CR embedding?
In this paper, among other things, we show how imposing the vacuum Einstein equations on an
associated space time can single out elements of this thin set. 

\subsection{The dual formulation}

Since we have $Z\dz\bar{Z}\neq 0$ we
can supplement these two complex vector fields by one
\emph{real} vector field $Z_0$ on $M$ so that the triple 
$(Z_0,Z,\bar{Z})$ constitutes a basis for complex vector fields on 
$M$. 
Associated to the basis $(Z_0,Z,\bar{Z})$ there is its dual basis of
1-forms $(\lambda,\mu,\bar{\mu})$ on $M$ satisfying 
$$Z_0\hook\lambda=1,\quad\quad Z\hook\mu=1,\quad\quad
\bar{Z}\hook\bar{\mu}=1,$$
all other contractions being identically equal to zero. Note that the
form $\lambda$ is real valued. 

The forms $(\lambda,\mu,\bar{\mu})$ are determined by the CR structure 
$\cs$ up to the following transformations:
\be
\lambda\mapsto\lambda'=f\lambda,\quad\quad\mu\mapsto\mu'=h\mu+p\lambda,\quad\quad\bar{\mu}\mapsto\bar{\mu}'=\bar{h}\bar{\mu}+\bar{p}\lambda,
\label{tra}
\ee
where $f\neq 0$ (real) and $h\neq 0, p$ (complex) are functions on
$M$.  

It is obvious that the differential equation
$$\der\zeta\dz\lambda\dz\mu=0$$
for a complex valued function $\zeta$ on $M$ is invariant under
the transformations (\ref{tra}). It is the dual version of the tangential
CR equation (\ref{hl}) and its solutions are just CR functions.

The Levi form of a 3-dimensional CR structure $\cs$ 
at a point is a nonvanishing multiple
of the value of the real valued function $\omega$ defined by 
\be
\lambda'\dz\der\lambda'=i\omega\lambda'\dz\mu'\dz\bar{\mu}'.\label{levi}
\ee
Thus the CR structure is strictly pseudoconvex iff 
$$\lambda\dz\der\lambda\neq 0$$
on $M$. (Note that this statement is independent of the choice of the
representative $\lambda'$ in (\ref{tra}).)

This enables us to formulate an equivalent definition of
a 3-dimensional CR structure, the original one, that was actually used
by Elie Cartan \cite{carcr}:
\begin{definition}
A 3-dimensional CR structure is a 3-dimensional real manifold $M$
equipped with an equivalence class of pairs of 1-forms $(\lambda,\mu)$
such that:
\begin{itemize}
\item $\lambda$ is real valued, $\mu$ is complex valued
\item $\lambda\dz\mu\dz\bar{\mu}\neq 0$ at each point of $M$
\item two pairs $(\lambda,\mu)$ and $(\lambda',\mu')$ are in the same
  class iff there exists functions $f$, $h$, $p$ such that (\ref{tra}) holds.
\end{itemize}
We will use this definition in the following and, rather than $(M,{\rm
  T}^{(1,0)})$, we will write
$(M,(\lambda,\mu))$ to stress that a CR structure is associated with a class
$[(\lambda,\mu)]$. 

\end{definition}  
\subsection{Lifting CR manifolds to Lorentzian spacetimes}\label{bub}
Three dimensional CR structures are very closely related to 
the so called congruences\footnote{In mathematical language the
  physicists' term `congruence' means: `foliation of a manifold by
  curves'; in this particular situation it means: `foliation by means
  of a three parameter family of null and shearfree geodesics'.} of null geodesics without shear in a
  spacetime \cite{PenRin,RobinsonTrautman,Taf}. These are well known in general relativity theory and
  proved to be very useful in the process of constructing nontrivial
  solutions to the vacuum Einstein equations in 4-dimensional
  manifolds equipped with Lorentzian metrics.      

Given a 3-dimensional CR manifold $M$, we consider a
representative $(\lambda,\mu)$ of the class  $[(\lambda,\mu)]$
defining the CR structure on it. On the 
Cartesian product ${\mathcal M}=M\times \bbR$ we have a distinguished
field of directions $k$, which is tangent to the $\bbR$ factor in
$\mathcal M$.  The four manifold ${\mathcal M}=M\times\bbR$ 
naturally projects onto $M$ with a projection $\pi:{\mathcal
M}\to M$ and $\pi_*(k)=0$. We choose a coordinate $r$ along the $\bbR$
factor, so that that $k$ may be represented by $k=\partial_r$. 
Omitting the pullbacks $\pi^*$ when expressing the forms on
$\mathcal M$, i.e. for example, denoting by $\mu$ the pullback
$\pi^*(\mu)$, 
we equip $\mathcal M$ with a class of metrics \cite{RobinsonTrautman,RoT,traa} 
\be
g=2P^2~[~\mu\bar{\mu}+\lambda(\der
r+W\mu+\bar{W}\bar{\mu}+H\lambda)~],\label{met}
\ee
where $P\neq 0, H$ (real) and $W$ (complex) are arbitrary functions on
$\mathcal M$; the expressions like e.g. $\mu\bar\mu$ denote the symmetrized tensor product:
$\mu\bar{\mu}=\tfrac{1}{2}(\mu\otimes\bar{\mu}+\bar{\mu}\otimes\mu)$.
We note that that the coordinate $r$ has no 
geometrical meaning; it can be replaced by any other function
$r'=r'(r,x,y,z)$ such that $\partial r'/\partial r\neq 0$.  

Now we consider the entire \emph{class} of metrics (\ref{met}), 
which depends on
arbitrary functions $P,W,H$ and the class of coordinates $r'$. 
We claim that this \emph{class} is \emph{naturally} attached to the CR
structure $(M,(\lambda,\mu))$. To see this, start with another
representative $(\lambda',\mu')$ of the class 
$[(\lambda,\mu)])$. The new forms
$(\lambda',\mu')$ are related to the previous choice
$(\lambda,\mu)$ via (\ref{tra}). Now, maintaining the same $P,W,H$
and $r$, write a metric $g'$ using formula (\ref{met}) with 
$(\lambda,\mu,\bar{\mu})$ replaced by $(\lambda',\mu',\bar{\mu}')$. 
Using definitions (\ref{tra}) the metric $g'$ can be reexpressed in terms
of the original forms $(\lambda,\mu,\bar{\mu})$. A short calculation
shows that after this, $g'$ has again form $(\ref{met})$ with merely
the functions $P,W,H$ and the coordinate $r$ being changed. Thus with each 
3-dimensional CR structure $(M,(\lambda,\mu))$ there is an associated  
4-dimensional manifold ${\mathcal M}=M\times \bbR$ with a 
class of Lorentzian metrics $g$ as in (\ref{met}). Now we can compare
$g$ and $g'$. It follows 
that there exists a nonvanishing real
function $\alpha$ and a 1-form $\varphi$ on $\mathcal M$ such that   
\be g'=\alpha^2 g+2 g(k)\varphi.\label{eq}
\ee 
Here $g(k)$ is the 1-form on $\mathcal M$ such that $X\hook g(k)=g(k,X)$.  
The class of Lorentzian 4-dimensional metrics $[g]$ with the
equivalence relation $g\sim g'$ iff $g$ and $g'$ are related by
(\ref{met}), (\ref{eq}) is called the class of metrics \emph{adapted} to the CR
structure $(M,(\lambda,\mu))$.

In each of the metrics from 
this class the lines tangent to the integral curves of the vector field 
$k=\partial_r$, are \emph{null}. They have the further property of
satisfying 
\be 
{\mathcal L}_kg=\Theta g+2g(k)\vartheta,\label{sf}
\ee
with a real function $\Theta$, the \emph{expansion}, 
and a real 1-form $\vartheta$ on $\mathcal
M$. This equation in particular means that the integral curves of $k$
are \emph{geodesics}. It also implies that the curves are
\emph{shearfree}, meaning that they preserve the natural conformal metric
defined by the class $[g]$ in the quotient space $k^\perp/k$. 

In the traditional language of physicsts ${\mathcal M}=M\times \bbR$ is equipped with a \emph{congruence} of null and
  shearfree geodesics tangent to $k$. Physicists say that this
  congruence is \emph{diverging} at a point of $\mathcal M$ 
iff the expansion $\Theta$ in (\ref{sf}) is nonvanishing at this point.   

The leaf space of integral curves of the congruence 
generated by $k$ can be identified with $M$. The property of the
  congruence of being null geodesic and shearfree means precisely that the
  3-dimensional leaf space $M$ of its integral curves has an abstract 
CR structure.

The above described procedure
of associating a metric $g$ from $[g]$, to a 
3-dimensional CR structure $(M,(\lambda,\mu))$, will be called a 
\emph{lift of the CR structure to a spacetime} 
\cite{rot,RobinsonTrautman,traa}.    

Given a lift of a CR structure $(M,(\lambda,\mu))$ to a spacetime, we
now briefly define two concepts, which are important for the
formulation of our main results. More detailed definitions are to be
found in Section \ref{sedf} and Section \ref{maxz}, respectively.

The first notion (see the end of Section \ref{sedf}) is defined as follows:

Having chosen a representative
$(\lambda,\mu)$ of $[(\lambda,\mu)]$ we pull it back to ${\mathcal
  M}=\bbR\times M$ by
$\pi^*$. Then we observe that the 2-dimensional complex distribution $N$, 
consisting of complex vector fields $Y$ on $\mathcal M$ satisfying $Y\hook
\pi^*(\lambda\dz\mu)=0$ 
is well defined. This is
because the 2-forms $\pi^*(\lambda\dz\mu)$ and $\pi^*(\lambda'\dz\mu')$
corresponding to different representatives of the class
$[(\lambda,\mu)]$ merely differ by the scale of a nonvanishing complex function.
The
distribution $N$ is called the \emph{distribution of} $\alpha$
\emph{planes} associated with the CR structure $(M,(\lambda,\mu))$ 
in the lifted spacetime $\mathcal M$.

The second notion we define here (again skipping the details to
Section \ref{maxz}, which includes a beatiful example due to Ivor Robinson) is a \emph{null Maxwell field aligned with the
  congruence}. Sometimes we will use the term `a \emph{null aligned
  Maxwell field}', for short. Suppose that in the class of forms
$[(\lambda,\mu)]$ there is a pair
$(\lambda',\mu')$ with the property that $\der(\pi^*(\lambda'\dz\mu'))=0$. 
Then we say that the CR structure $(M,(\lambda,\mu))$ admits a null 
aligned Maxwell field ${\mathcal F}=\pi^*(\lambda'\dz\mu')$.  

\subsection{From congruences of shearfree and null geodesics to CR manifolds}\label{frc}
The fact that any 3-dimensional CR manifold $(M,(\la,\mu))$ locally 
defines a class of Lorentzian metrics (\ref{met}) on $M\times\bbR$ 
in which the lines tangent to the $\bbR$ factor are null and shearfree geodesics has also its converse. This is given by the following theorem \cite{PenRin,rot,RobinsonTrautman,so,Taf,traa}.
\begin{theorem}\label{rttr}
Let $({\mathcal M},g)$ be a 4-dimensional manifold equipped with a Lorentzian metric. Suppose that $\mathcal M$ is foliated by a 3-parameter family of curves which are null geodesics without shear. Then $\mathcal M$ is locally a cartesian product ${\mathcal M}=M\times\bbR$ with $M$ being a 3-dimensional CR manifold. The CR 
structure $(M,(\lambda,\mu))$ on $M$ is uniquely determined by $({\mathcal M},g)$ and the shearfree congruence on $\mathcal M$. If $r$ is a real coordinate 
such that $k=\partial_r$ is tangent to the congruence, then the 
Lorentzian metric $g$ on $\mathcal M$ can be locally represented by 
(\ref{met}) with some specific functions $P,W,H$ depending on the choice of the representatives $(\lambda,\mu)$ of the corresponding CR structure. 
\end{theorem}     

It should be noticed that 4-dimensional Lorentzian metrics of the form (\ref{met}) have been studied by physicists since the late 1950s. In 
particular physicists found a lot of examples of such metrics which 
satisfy Einstein's equations $Ric(g)=\Lambda g$. Among them is the 
Schwarzschild solution (the corresponding CR manifold is Levi flat), the  
Taub-NUT solution (the corresponding CR manifold is almost everywhere 
CR equivalent to the Heisenberg group CR structure), the Kerr rotating 
black hole solution (the corresponding CR manifold is almost
everywhere  
strictly pseudoconvex and has only a 2-dimensional group of local CR symmetries; so it is not CR equivalent to the Heisenberg group CR structure), and many others (see \cite{KSMH}, sections devoted to algebraically 
special solutions, for huge 
families of Einstein examples with various CR structures).  Physicists 
were looking for the Einstein solutions among the metrics (\ref{met}) 
because such metrics were believed to correspond to gravitational radiation. 
Although the understanding of the mathematical fact that there is a CR structure behind the scenes came much later, physicists were from the very begining 
aware that radiative Maxwell or Einstein 
fields impose a sort of a complex structure in the underlying
spacetimes. The notion of a CR structure was implicit in such papers as  \cite{Debney,GS,Kerr,Newman,ps,Penrose,rob,robin,RT,RR,Sa,Sachs,so,Sommers,Tr8,Tr2}, but physicists  
did not manage to abstract the concept for about twenty years \cite{Sommers,Taf,traa}. 
Ironically, at about the same time when the systematic work on
gravitational radiation started,
the notion of a CR structure was being revived in mathematics, 
due to the discovery of
the nonsolvability of equations of Hans Lewy type \cite{Le2}. 
Mathematicians were
however unaware of the developement in general relativity theory and also did
not make the connection. We hope that this paper fills the gap between
these two areas of mathematics and physics.

Our main motivation is the paper \cite{lnt} and the research on relations
between 3-dimensional CR structures and algebraically special
gravitational fields undertaken by the Warsaw Relativity Group in the
1980s
\cite{lew,lewPhD,ln,ln1,lnt1,nuro,nurta,nurphd,RobinsonTrautman,rot,RoT,Taf,Taff,traa,Tr}. We are also
inspired and impressed by the works of relativists on
gravitational radiation; in particular by the contributions of Andrzej Trautman, Ivor Robinson, Roger
Penrose, Ray Sachs, Roy Kerr, 
Ted Newman and Jacek Tafel. 
\section{Local CR embeddability theorems}
In the sequel we assume our CR structures have a sufficiently high
\emph{finite} order of differentiability, in particular they need \emph{not} be real
analytic. All considerations are \emph{local} on $M$.
\begin{theorem}\label{firstcr}
Let $M$ be a sufficiently smooth 3-dimensional CR manifold. If it has a  
lift to a spacetime whose congruence of null and shearfree
geodesics is diverging over the points where the Levi form vanishes, and
whose complexified Ricci tensor vanishes on the 
distribution of $\alpha$ planes associated with the congruence,   
then $M$ admits a CR function $\zeta$ such that
$\der\zeta\dz\der\bar{\zeta}$ is nowhere zero.
\end{theorem}
In particular Theorem \ref{firstcr} applies to the strictly
pseudoconvex case. Thus it rules out the generic situation, like the example of
Nirenberg, in which all CR functions are locally constant. 

Actually, in the strictly pseudoconvex case, we have a stronger theorem. 
\begin{theorem}\label{fistcr}
Let $M$ be a sufficiently smooth 
strictly pseudoconvex 3-dimensional CR manifold. Then the following
two conditions are equivalent.
\begin{itemize}
\item[(i)] $M$ admits a lift to a spacetime whose complexified Ricci
  tensor vanishes on the corresponding distribution of $\alpha$ planes
\item[(ii)] $M$ admits one CR function $\zeta$ such that
  $\der\zeta\dz\der\bar{\zeta}$ is nowhere zero.
\end{itemize}
\end{theorem} 
Our next result gives an if and only if criterion for the
\emph{local embeddability} of a sufficiently smooth, \emph{not} necessarily
real analytic, strictly pseudoconvex 3-dimensional CR manifold: 
\begin{theorem}\label{maxc}
Let $M$ be a sufficiently smooth strictly pseudoconvex 3-dimensional CR manifold.
It is locally CR embeddable as a hypersurface in $\bbC^2$ if and only if 
\begin{itemize}
\item[(i)] it  
admits a lift to a spacetime
whose complexified Ricci tensor vanishes on the corresponding
distribution of $\alpha$ planes, 
\item[(ii)] and it admits a nontrivial null
Maxwell field aligned with the null congruence of shearfree geodesics
corresponding to the CR structure on $M$. 
\end{itemize}
\end{theorem} 
Next we abandon the requirement about the existence of a null 
Maxwell field and replace it
by further assumptions about the curvature of the lifted
spacetime. This leads to a remarkable theorem. 
\begin{theorem}\label{petrovt}
Let $M$ be a sufficiently smooth 
strictly pseudoconvex 3-dimensional CR manifold. If $M$ admits a lift
to a Ricci flat spacetime 
which has locally
constant 
Petrov
type, then it is locally CR embeddable as a hypersurface in $\bbC^2$.
\end{theorem}
The hypothesis of \emph{locally constant Petrov type} is a technical
assumption, which will be explained in detail in Section
\ref{sedf}. Here we only mention that according
to the theory of exact solutions of Einstein equations, a spacetime
metric at a point can have one of the six Petrov types \cite{caru,Petrov,ps}. These are:
Petrov types $I$, $II$, $D$, $III$, $N$ or $0$, and they may change
from point to point in the spacetime. With this information
we may formulate the following theorem.  
\begin{theorem}\label{petovt}
Let $(M,(\lambda,\mu))$ be a real analytic strictly pseudoconvex 
3-dimensional CR manifold. Then it always has a lift to a spacetime
satsifying the Einstein equations $Ric(g)=\Lambda g+ \Phi
\lambda\otimes\lambda$ whose Petrov type is $II$ or $D$.  
\end{theorem}
Here $\Lambda$ is the \emph{cosmological constant} and the function $\Phi$
corresponds to the energy momentum tensor of \emph{pure radiation}. We
believe that this theorem is also true when we replace the term `real
analytic' with `sufficiently smooth embeddable' (see Remark \ref{rekm}).\\

The proofs of the above theorems will emerge in the discussion which follows.
Actually the theorems stated above, constitute only a selection of the
results we prove in the paper. In special cases, which are
systematically studied in the main body of the paper, we obtain
sharper results than stated here.\\

We close this section with a remark about the nontriviality of 
Theorems \ref{firstcr} and \ref{petrovt}. As was already mentioned 
after Theorem \ref{rttr} both 
theorems are far from being empty. There is an abundance of Ricci flat 
Lorentzian 4-metrics which admit a congruence of null and shearfree 
geodesics. The encyclopedia book \cite{KSMH} gives an up to date catalog 
of such metrics in the sections devoted to algebraically special 
vacuum solutions (Sections 26 through 30 in the second edition). Every Ricci 
flat metric in these sections of the book has a corresponding 
3-dimensional CR structure. This may be Levi flat everywhere 
(as is the case for the Schwarzschild metric) or strictly pseudoconvex in 3-dimensional regions and Levi flat on some lower dimmensional sets as in the following example:

We consider the metric  
$$
g=2\Big({\mathcal P}^2\mu\bar{\mu}+\lambda(\der r'+{\mathcal
  W}\mu+\bar{\mathcal W}\mu+{\mathcal H}\lambda)\Big),$$
where  
$$\la=\der
u+\frac{i\big(2b+(a+b)\zeta\bar{\zeta}\big)}{\zeta(1+\tfrac{\zeta\bar{\zeta}}{2})^2}\der
\zeta-\frac{i\big(2b+(a+b)\zeta\bar{\zeta}\big)}{\bar{\zeta}(1+\tfrac{\zeta\bar{\zeta}}{2})^2}\der
\bar{\zeta},\quad\quad\mu=\der \zeta,$$
$${\mathcal P}^2=\frac{r'^2}{(1+\tfrac{\zeta\bar{\zeta}}{2})^2}+\frac{\big(b-a+(b+a)\tfrac{\zeta\bar{\zeta}}{2}\big)^2}{(1+\tfrac{\zeta\bar{\zeta}}{2})^4},$$
$$
{\mathcal W}=\frac{ia\bar{\zeta}}{(1+\tfrac{\zeta\bar{\zeta}}{2})^2},\quad\quad
{\mathcal H}=-\frac{1}{2}+\frac{mr'+b^2-ab~\frac{1-\tfrac{\zeta\bar{\zeta}}{2}}{1+\tfrac{\zeta\bar{\zeta}}{2}}}{r'^2+\frac{\big(
  b-a+(b+a)\tfrac{\zeta\bar{\zeta}}{2}\big)^2}{(1+\tfrac{\zeta\bar{\zeta}}{2})^2}},$$
and $m,a,b$ are \emph{real} constants.

Clearly this metric is in the 
form (\ref{met}), and as such may be considered as the lift of a 
CR structure $(\lambda,\mu)$ which is defined on 
the 3-dimensional manifold $M$ parametrized by 
$(u,{\rm Re}(\zeta),{\rm Im}(\zeta))$. 

The interesting feature of this scarry-looking 3-parameter family of 
metrics is that it is \emph{Ricci flat} for all values of the coordinates $(u,{\rm Re}(\zeta),{\rm Im}(\zeta),r')$ in which $g$ is not singular \cite{KSMH}. 
Actually if $b=0$ the above metric is just the Kerr rotating black hole with mass 
$m$ and the angular momentum parameter $a$; if $a=b=0$ the metric $g$ describes the Schwarzsschild black hole with mass $m$. If $m=a=0$ the corresponding metric is the Taub-NUT vacuum metric. 

Calculating $\lambda\dz\der\lambda$ we get:
$$\la\dz\der\la=i\frac{(a+b)\zeta\bar{\zeta}-2(a-b)}{(1+\tfrac{\zeta\bar{\zeta}}{2})^3}\der u\dz\der\zeta\dz\der\bar{\zeta}.$$

This means that for each value of the three parameters $(m,a,b)$ 
the corresponding CR structure $(M,(\la,\mu))$ is strictly 
pseudoconvex everywhere except the points where 
$$(a+b)\zeta\bar{\zeta}-2(a-b)=0.$$

Note that if $a>b$ and $a\neq -b$, there is an entire cylinder 
$\zeta\bar{\zeta}=\frac{2(a-b)}{a+b}$ in $M$ on which 
$\lambda\dz\der\lambda=0$. In such case the corresponding CR structure is Levi flat on this cylinder and strictly pseudoconvex outside it. 
A short calculation shows that on this cylinder 
the shearfree congruence of null geodesics tangent to
$k=\partial_{r'}$ is diverging everywhere. So this case is a 
nontrivial example of a metric which appears in Theorem \ref{firstcr}. Many more such examples can be found in the appropriate sections of \cite{KSMH}.
     
\section{The Einstein equations and CR functions}\label{sefd}
\subsection{The first CR function}\label{sedf}
Here we prove Theorem \ref{firstcr} by adapting  
the argument presented in \cite{KSMH}. 

We consider a general 4-dimensional spacetime $\mathcal M$, i.e. a 4-manifold
equipped with a metric $g$ of Lorentzian signature $(+,+,+,-)$. We assume
that the spacetime $\mathcal M$ admits a null congruence of
shearfree and null geodesics. This may be described as
follows:

The congruence is tangent to a nonvanishing vector field $k$ which is
\emph{null}, $g(k,k)=0$. Having $k$, we introduce a
coframe $(\theta^1,\theta^2,\theta^3,\theta^4)$ on $\mathcal M$ such
that 
\be
g=2(\theta^1\theta^2+\theta^3\theta^4),\quad\quad
\theta^3=g(k,\cdot),\quad\quad{\rm
  and}\quad\quad k\hook\theta^1=k\hook\theta^2=0.
\label{cof}
\ee
Note that, due to the signature of the metric $g$, this definition
implies that the 1-forms $\theta^3$, $\theta^4$ are \emph{real} valued, whereas
the 1-forms $\theta^1$, $\theta^2$ are \emph{complex} valued with 
$\bar{\theta}^2=\theta^1$. Note also that the coframe 
$(\theta^1,\theta^2,\theta^3,\theta^4)$ is \emph{not} uniquely 
defined by (\ref{cof}). It is given up to a linear transformation
associated with a 4-dimensional parabolic 
subgroup of the Lorentz group preserving the null direction $k$. Explicitly:
\begin{eqnarray}
&&{\theta^1}\phantom{}'={\rm e}^{i\phi}(\theta^1+\bar{B}\theta^3)\nonumber\\
&&{\theta^2}\phantom{}'={\rm
  e}^{-i\phi}(\theta^2+B\theta^3)\label{cofr}\\
&&{\theta^3}\phantom{}'=A\theta^3\nonumber\\
&&{\theta^4}\phantom{}'=A^{-1}(\theta^4-B\theta^1-\bar{B}\theta^2-B\bar{B}\theta^3),\nonumber
\end{eqnarray}
where $A\neq 0$, $\phi$ are real functions, and $B$ is a complex
function. The coframes (\ref{cofr}) are said to be \emph{adapted} to $k$. 

Imagine now that we have $k$, which in some adapted
coframe $(\theta^1,\theta^2,\theta^3,\theta^4)$ satisfies:
\begin{eqnarray}
&&\der\theta^3\dz\theta^1\dz\theta^3=-\bar{\kappa}~\theta^1\dz\theta^2\dz\theta^3\dz\theta^4\label{kass}
\\
&&\der\theta^1\dz\theta^1\dz\theta^3=-\bar{\sigma}~\theta^1\dz\theta^2\dz\theta^3\dz\theta^4,\nonumber
\end{eqnarray}
with some complex functions $\kappa$ and $\sigma$. Then, it follows 
that in any other adapted coframe (\ref{cofr}) we have
\begin{eqnarray*}
&&\der{\theta^3}\phantom{}'\dz\theta^1\phantom{}'\dz\theta^3\phantom{}'=-A^2{\rm
  e}^{i\phi}\bar{\kappa}~\theta^1\phantom{}'\dz\theta^2\phantom{}'\dz\theta^3\phantom{}'\dz\theta^4\phantom{}'
\\
&&\der\theta^1\phantom{}'\dz\theta^1\phantom{}'\dz\theta^3\phantom{}'=-A{\rm
  e}^{2i\phi}(\bar{B}\bar{\kappa}+\bar{\sigma})~\theta^1\phantom{}'\dz\theta^2\phantom{}'\dz\theta^3\phantom{}'\dz\theta^4\phantom{}',
\end{eqnarray*}
as can be easily checked. This means that the \emph{simultaneous vanishing}
or \emph{not} of both coefficients $\kappa$, $\sigma$ at a point, is independent of the
choice of the adapted coframe, and thus is a property of a null
congruence associated with $k$. If 
\be
\kappa=\sigma=0\label{kasi}
\ee
everywhere, the null congruence associated with $k$ is called a \emph{shearfree} congruence of null \emph{geodesics}. We mention in
addition that the vanishing of $\kappa$ alone at a point, is also an invariant
property of the congruence. If we assume nothing about $\sigma$ but
require that $\kappa=0$ everywhere, such
a congruence turns out to consist of null \emph{geodesics}, which may
or may not be shearfree.

It is also worthwhile to look at the transformations of
$\der\theta^3\dz\theta^3$. To be consistent with (\ref{kass}) we write
it as:
\be
\der\theta^3\dz\theta^3=i\Omega\theta^1\dz\theta^2\dz\theta^3-(\kappa\theta^1+\bar{\kappa}\theta^2)\dz\theta^3\dz\theta^4,\label{twii}
\ee
with a real function $\Omega$ on $\mathcal M$. Changing the adapted
coframe to (\ref{cofr}) we get:
$$
\der\theta^3\phantom{}'\dz\theta^3\phantom{}'=
A(\bar{B}\kappa-B\bar{\kappa}+i\Omega)\theta^1\phantom{}'\dz\theta^2\phantom{}'\dz\theta^3\phantom{}'-A^2\big({\rm
  e}^{-i\phi}\kappa\theta^1\phantom{}'+{\rm
  e}^{i\phi}\bar{\kappa}\theta^2\phantom{}')\dz\theta^3\phantom{}'\dz\theta^4\phantom{}'.
$$
In the case of a geodesic null congruence, $\kappa\equiv 0$, these equations reduce
respectively to:
\be
\der\theta^3\dz\theta^3=i\Omega\theta^1\dz\theta^2\dz\theta^3\label{twi}
\ee 
and 
$$\der\theta^3\phantom{}'\dz\theta^3\phantom{}'=iA\Omega\theta^1\phantom{}'\dz\theta^2\phantom{}'\dz\theta^3\phantom{}'.$$
This proves that in such case the vanishing or not of
$\Omega$ at a point is
also an invariant property of $k$, and thus, of the
congruence. Consider a \emph{null} and 
\emph{geodesic} congruence at a point $x\in{\mathcal M}$. If 
$$\Omega\neq 0$$
at $x$, we say that the congruence is \emph{twisting} there. If 
$$\Omega=0$$
at $x$, we say that the
congruence is \emph{not twisting} at $x$.

From now on we assume that we have $k$ with $\kappa=\sigma=0$ everywhere,
i.e. that we have a null congruence of shearfree geodesics in
$\mathcal M$. Our next step is to give the geometric interpretation of
(\ref{kasi}).  

Choosing an adapted coframe (\ref{cof}), using Cartan's formula 
${\mathcal L}_k\theta=k\hook\der\theta+\der(k\hook\theta)$
 for the Lie derivative of a 1-form $\theta$, and the respective 
equations (\ref{twi}), (\ref{kass})$_2$ with $\sigma=0$, we easily get
$$({\mathcal L}_k\theta^3)\dz\theta^3=0,\quad\quad{\rm and}\quad\quad
({\mathcal L}_k\theta^1)\dz\theta^1\dz\theta^3=0,$$
everywhere in $\mathcal M$. The meaning of these two equations is
obvious: the real 1-form $\theta^3$, when Lie transported by 
the flow $\phi_t$ generated by the congruence, transforms as
$\phi_t^*(\theta^3)=f\theta^3$, where $f$ is a real function on
$I\times{\mathcal M}$, $t\in I$; the complex 1-form $\theta^1$ transforms as
$\phi_t^*(\theta^1)=h\theta^1+p\theta^3$, where $h$, $p$ are complex
functions on $I\times{\mathcal M}$. Since $\phi_t$ is a local diffeomorphism,
the functions $f$ and $h$ are locally nonvanishing.

Now, taking any hypersurface $M$ in $\mathcal M$ transversal to $k$, we
equip it with a CR structure $(\lambda',\mu',\bar{\mu}')$ as in
(\ref{tra}) by setting 
$$\lambda'=\iota^*(\theta^3),\quad
\mu'=\iota^*(\theta^1),\quad
\bar{\mu}'=\iota^*(\theta^2).$$ Here $\iota:M\to\mathcal
M$ is the natural inclusion of $M$ in $\mathcal M$, so that
$\iota^*\theta$ is just the restriction of $\theta$ to $M$. The above discussed
changes of $\theta^3$ and $\theta^1$ along $k$ imply that the CR
structures on any two transversal hypersurfaces are CR
equivalent. The pseudoconvexity property of these CR structures is
easily described by means of equation (\ref{twi}). Indeed, pulling back
this equation by means of $\iota^*$, from the spacetime to
$M$, 
we get
$$\der\lambda'\dz\lambda'=i\iota^*(\Omega)\mu'\dz\bar{\mu}'\dz\lambda'.$$
This means that the Levi form $\omega$, as defined by (\ref{levi}), is
$\omega=\iota^*(\Omega)$. Thus the CR structure $(\lambda',\mu',\bar{\mu}')$
has a nonvanishing Levi form $\omega\neq 0$ at $p\in M$ iff  
the unique congruence curve passing through $p$ is twisting at $\iota(p)$.     

Summarizing we have the following lemma \cite{RobinsonTrautman}. 
\begin{lemma}
The 3-dimensional leaf space of a null congruence of shearfree
geodesics in spacetime is locally a CR manifold. This CR manifold is
strictly pseudoconvex at a point iff the congruence curve is
twisting at the corresponding point in spacetime.
\end{lemma}   

It is convenient to choose a real function $r$ on $\mathcal M$, so
that $k=P^{-1}\partial_r$, where  $P\neq 0$ is a real function on
$\mathcal M$. Then our Lemma guarantees that the adapted
coframe $(\theta^1,\theta^2,\theta^3,\theta^4)$ can be chosen in such
a way that
\begin{eqnarray}
&&\theta^1=P\mu\nonumber\\
&&\theta^2=P\bar{\mu}\nonumber\\
&&\theta^3=P\lambda\label{choco}\\
&&\theta^4=P(\der r+W\mu+\bar{W}\bar{\mu}+H\lambda),\nonumber
\end{eqnarray}  
where locally ${\mathcal M}=\bbR\times M$, $\lambda,\mu$ are the 1-forms
on $M$ which define the CR structure there, and the functions
$H$ (real) and $W$ (complex) are functions on ${\mathcal M}$. This in
particular means that in addition to $k\hook\lambda=k\hook\mu=0$ we
also have $k\hook\der\lambda=k\hook\der\mu=0$. Thus we just
demonstrated how a shearfree congruence of null geodesics in spacetime
restricts the spacetime metric to the form (\ref{met}). As we already
mentioned in Section \ref{bub}, we also have the statement in the opposite
direction: given a 3-dimensional CR structure represented on $M$ via
forms $(\lambda,\mu,\bar{\mu})$ we lift it to a spacetime ${\mathcal
  M}=\bbR\times M$ with metric (\ref{met}) and with a null congruence
of shearfree geodesics represented by $k=P^{-1}\partial_r$.  

Whether we start with a 3-dimensional CR structure and then define the
spacetime with a null congruence of shearfree geodesics, or
immediately start with a spacetime with such congruence, we may try to
impose some curvature conditions on $g$. The question arises if these
conditions say something about the underlying CR geometry. 

To study this question we consider Cartan's structure equations for the metric (\ref{met}),
written in an adapted coframe (\ref{choco}). These are:
\begin{eqnarray}
&&\der\theta^i+\Gamma^i_{~j}\dz\theta^k=0\label{ca1}\\
&&\der\Gamma^i_{~j}+\Gamma^i_{~k}\dz\Gamma^k_{~j}=\tfrac12 R^i_{~jkl}\theta^k\dz\theta^l.\label{ca2}
\end{eqnarray}
Here the 1-forms $\Gamma^i_{~j}$ 
are the Levi-Civita connection 1-forms. They define 
$\Gamma_{ij}=g_{ik}\Gamma^k_{~j}$, which 
satisfy $\Gamma_{ij}=-\Gamma_{ji}$. Modulo the symmetry the only
nonzero components of the metric are $g_{12}=g_{34}=1$. Its inverse 
$g^{ij}$, again modulo symmetry, has $g^{12}=g^{34}=1$ as the only nonvanishing components. The coeffcients $R^i_{~jkl}$
  are the Riemann tensor coeffcients. The Ricci tensor is defined as
  $R_{ij}=R^k_{~ikj}$. Using the metric $g_{ij}$ we also
  define $R_{ijkl}=g_{im}R^m_{~jkl}$. This can be used to define the
  covariant components of the Weyl tensor $C^i_{~jkl}$ via:
\be
C_{ijkl}=R_{ijkl}+\tfrac16 R
  (g_{ik}g_{lj}-g_{il}g_{kj})+\tfrac12(g_{il}R_{kj}-g_{ik}R_{lj}+g_{jk}R_{li}-g_{jl}R_{ki}).\label{wyel}\ee
Here $R=R_{ij}g^{ij}$ is the Ricci scalar. The Weyl tensor $C^i_{jkl}=g^{im}C_{mjkl}$ carries the conformal
  information about the spacetime.

\begin{remark}\label{rek}
Note that, due to the reality conditions: $\bar{\theta}^1=\theta^2$,
$\bar{\theta}^3=\theta^3$, $\bar{\theta}^4=\theta^4$ we have
$\bar{R}_{13}=R_{23}$, $\bar{R}_{33}=R_{33}$,
$\bar{\Gamma}_{14}=\Gamma_{24}$, etc. This means that 
complex conjugation interchanges the 
indices $1\leftrightarrow 2$ and leaves the indices $3$ and $4$
unchanged.
\end{remark}

Now, since $\kappa=\sigma=0$, the connection 1-form $\Gamma_{24}$ satisfying Cartan's
first structure equation (\ref{ca1}) and, being compatible with our
conventions (\ref{kass}) and (\ref{twii}), must be a linear combination
of $\theta^1$ and $\theta^3$ only:
\be
\Gamma_{24}=-\bar{\rho}\theta^1-\bar{\tau}\theta^3.\label{wyco}
\ee
This equation defines complex functions $\rho$ and $\tau$. Our
conventions do not give any restrictions on $\tau$. On the other hand,
to be compatible with (\ref{sf}) and (\ref{twi}) the function $\rho$ has
to assume the form  
$$\rho=\tfrac12(-\Theta+i\Omega),$$
with real $\Omega$ being the twist, and real $\Theta$ being 
the expansion of the congruence. 

  
Now we pass to the curvature analysis. The first option is to 
impose the Einstein equations 
\be
R_{ij}=\Lambda g_{ij},\quad\quad i,j=1,2,3,4\label{eie}
\ee
on $g$. These 
equations conveniently split into three types of equations \cite{KSMH}:
\begin{itemize}
\item[(a)] $R_{22}=R_{24}=R_{44}=0$, 
\item[(b)] $R_{12}=R_{34}(=\Lambda)$,
\item[(c)] $R_{33}=R_{23}=0$.
\end{itemize}   
According to Remark \ref{rek} the set (a) consists of five real equations
($R_{44}$ is real!), the set
(b) consists of two real
equations and the set (c) consists of three real equations.\\

To proceed with the proof of Theorem \ref{firstcr} we now 
focus our attention on equations (a).\\

First, using the definition of the Ricci tensor, we write equations (a)
in terms of the Riemann tensor coefficients. A short calculation shows
that, modulo the symmetries of the Riemann tensor, they are equivalent to:
\begin{eqnarray*}
&&R_{44}=0\quad\Leftrightarrow\quad R_{2414}=0\\
&&R_{24}=0\quad\Leftrightarrow\quad R_{2412}-R_{2434}=0\\
&&R_{22}=0\quad\Leftrightarrow\quad R_{2423}=0.
\end{eqnarray*}
Second, we invoke the celebrated theorem due to Goldberg and Sachs
\cite{GS}. We will need this theorem in a very technical
version. Before giving its formulation suitable for our purposes it is
worthwhile to present its original meaning. 

As noted by Cartan \cite{caru}, every 4-dimensional spacetime
$\mathcal M$ distinguishes 
at each point at most \emph{four} null directions - the \emph{principal
null directions}, as they are nowadays called. If at a point
$p\in{\mathcal M}$ the Weyl tensor of the metric is not zero, then the
number $s$ of distinct directions is $1\leq s\leq 4$. The number $s$ at a
point depends on the Weyl tensor in an algebraic fashion: 
it is the number of distinct roots of a certain fourth order complex 
polynomial associated with the Weyl tensor. If $s=4$ at
$p\in\mathcal M$ the spacetime is called \emph{algebraically general}
at $p$. If $1\leq s\leq 3$, the spacetime is called \emph{algebraically
special} at $p$. If $1\leq s\leq 3$, then at least two of Cartan's principal
null directions coincide. The coinciding principal null directions are
called \emph{multiple} principal null directions.

All the possible degeneracies of the principal
null directions at a point can be enumerated by the possible
partitions of the number \emph{four}. Thus we have cases $[1111]$, $[112]$,
$[22]$, $[13]$, $[4]$. 
The algebraically general case corresponds to $[1111]$. The case where
there is only one doubly
degenerate principal null direction corresponds to $[112]$.  The case with 
two doubly degenerate principal null direction coresponds to $[22]$, and so
on. This classification of points in a spacetime, is known as the
\emph{Petrov types} \cite{Petrov}, although it was known to Elie Cartan
\cite{caru}, and was brought to its contemporary form by Roger Penrose
\cite{ps}. In Penrose's formulation it reads as follows:
\vspace{.5cm}

\parbox{.5\linewidth}{
\hspace{-.6cm}\, \,\,(i) Petrov type I (`non-degenerate') \({\pmb [}\ts 1111\ts{\pmb ]}\),

\hspace{-.6cm}\, (ii) Petrov type II \({\pmb [}\ts 112\ts{\pmb ]}\),

\hspace{-.67cm}\, (iii) Petrov type III \({\pmb [}\ts 13\ts{\pmb ]}\),

\hspace{-.6cm} (iv) Petrov type D (`degenerate') \({\pmb [}\ts 22\ts{\pmb ]}\),

\hspace{-.5cm}\,\,(v) Petrov type N (`null') \( {\pmb [}\ts 4\ts{\pmb ]}\).}

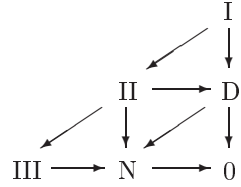
\parbox{.5\linewidth}{\begin{picture}(10,10)(-160,15)

\put(150,85){I}\put(152,82){\vector(0,-1){15}}
\put(110,55){II}\put(113,52){\vector(0,-1){15}}
\put(149,55){D}\put(152,52){\vector(0,-1){15}}
\put(70,25){III}\put(110,25){N}\put(150,25){0}
\put(85,29){\vector(1,0){20}}
\put(123,29){\vector(1,0){22}}
\put(123,59){\vector(1,0){22}}
\put(144,82){\vector(-3,-2){23}}
\put(104,52){\vector(-3,-2){23}}
\put(143,52){\vector(-3,-2){23}}
\end{picture}
}

\smallskip
\noindent The \(0\) in the diagram above, the \emph{Penrose diagram}
as it is called, represents a vanishing Weyl tensor at a point.
The arrows point towards more special cases. 

A convenient way to determine the Petrov type is to
calculate the \emph{Weyl scalars} $\Psi_0$, $\Psi_1$, $\Psi_2$, $\Psi_3$,
$\Psi_4$ at a point. These quantities are complex numbers at a point, which fully determine the 10 independent components of the Weyl tensor at this point. In a frame
$(e_1,e_2,e_3,e_4)=(m,\bar{m},l,k)$, dual to the coframe adapted to a 
null vector field $k$, they are defined by
\begin{eqnarray*}
&&\Psi_0=C_{ijkl}k^im^jk^km^l=C_{4141}=R_{4141}\\
&&\Psi_1=C_{ijkl}k^il^jk^km^l=C_{4341}=\tfrac12(R_{4341}+R_{1421})\\
&&\Psi_2=C_{ijkl}k^im^jl^k\bar{m}^l=C_{4132}\\
&&\Psi_3=C_{ijkl}l^ik^jl^k\bar{m}^l=C_{3432}=\tfrac12(R_{3432}+R_{2312})\\
&&\Psi_4=C_{ijkl}l^i\bar{m}^jl^k\bar{m}^l=C_{3232}=R_{3232}.
\end{eqnarray*}
For the convenience of the reader, we we have used here formula
(\ref{wyel}), to reexpress $\Psi_0$, $\Psi_1$, $\Psi_3$ and $\Psi_4$
in terms of the Riemann tensor components.

The importance of the Weyl scalars for the Petrov classification 
consists in the following observations:
\begin{itemize}
\item $\Psi_0=0\Longleftrightarrow k~ {\rm is~ a~ principal~ null~
  direction}$ 
\item $\Psi_0=\Psi_1=0,\quad \Psi_2\neq 0 \Longleftrightarrow k~
  {\rm has ~degeneracy}~ [112] ~{\rm or}~[22]$
\item $\Psi_0=\Psi_1=\Psi_2=0,\quad \Psi_3\neq 0 \Longleftrightarrow k~
  {\rm has ~degeneracy}~ [13]$
\item $\Psi_0=\Psi_1=\Psi_2=\Psi_3=0,\quad \Psi_4\neq 0
  \Longleftrightarrow 
k~{\rm has ~degeneracy}~ [4]$
\item $\Psi_0=\Psi_1=\Psi_2=\Psi_3=\Psi_4=0 \Longleftrightarrow {\rm
  Weyl~tensor~is ~zero.}$
\end{itemize}
We will use them in Sections \ref{ree} and \ref{reee}. 

In turns out, and
this has been known for years by physicists, that if the spacetime admits
a congruence of shearfree and null geodesics, then tangent
vectors to the congruence at a point should be aligned with one of the
principal null directions. The Goldberg-Sachs theorem says more \cite{GS}:
\begin{theorem}\label{ggss}(Goldberg-Sachs): Suppose that a 4-dimensional
  spacetime satisfies Einstein's equations $Ric(g)=\Lambda g$. Then
  the following conditions are equivalent:
\begin{itemize}
\item[i)] The spacetime admits a null congruence of shearfree geodesics
  tangent to a vector field $k$.
\item[ii)] The spacetime is algebraically special with a multiple principal
  null direction tangent to $k$.
\end{itemize}
\end{theorem}   
Algebraically special fields are important in physics because they
are connected with what \emph{gravitational radiation}
could be. By this we mean the following. 

Roughly speaking, if one observes a \emph{general} gravitational field \emph{far from the
sources} and measures the `distance' from the sources by means of
$r>0$, then studying the $r$-dependence of the Weyl tensor $C$ of
the metric (which in the
\emph{empty} spacetime describes the
gravitational field strength), he will discover a 
$$C=\frac{N}{r}+\frac{III}{r^2}+\frac{II+D}{r^3}+\frac{I}{r^4}+O(\frac{1}{r^5})$$
behavior, as $r\to\infty$.

Here $N$, $III$, $II$, $D$ and $I$ are tensorial quantities with all
the symmetries of the Weyl tensor. They have the respective 
algebraic Petrov type
denoted by the corresponding symbols. Thus, a \emph{general} 
gravitational field far
from the sources is of Petrov type $N$, as first observed in
\cite{Tr8}. Approaching the sources, the field becomes less and less algebraically special. This
is the so called \emph{peeling property} of the gravitational field 
\cite{Sa,Sachs}. It gives an \emph{algebraic} criterion for a
gravitational field to be `radiative'.\\    

Since the original paper of Goldberg and Sachs \cite{GS} the theorem
was strengthend in various ways. In particular, it is known that 
to achieve the implication i)$\Longrightarrow$ii) in Theorem \ref{ggss} the full set of the Einstein equations is not needed. One
can weaken $Ric(g)=\Lambda g$ to our equations (a) and the implication 
i)$\Longrightarrow$ii) in Theorem \ref{ggss} is still true. Also the
geometric condition about algebraic speciality can be reformulated
directly in terms of the vanishing of certain components of the Weyl tensor, and in turn, of the Riemann tensor. Such a form of the theorem is needed for proving our Theorem \ref{firstcr}. We quote it below \cite{ghn} (see also Lemma 2.2 on p. 577 in \cite{przp}):    
\begin{theorem}
Suppose that a 4-dimensional manifold $\mathcal M$ is equipped with a Lorentzian
metric which in a null coframe $(\theta^1,\theta^2,\theta^3,\theta^4)$
has the form
$$g=2(\theta^1\theta^2+\theta^3\theta^4),$$
with $\bar{\theta}^1=\theta^2$, $\bar{\theta}^3=\theta^3$, 
$\bar{\theta}^4=\theta^4$. Let $\kappa$ and $\sigma$ be given by
(\ref{kass}). 

Assume that the Ricci tensor of $g$, in the coframe 
$(\theta^1,\theta^2,\theta^3,\theta^4)$, satisfies  
$$R_{22}=R_{24}=R_{44}=0$$
everywhere on $\mathcal M$.
Then  
\begin{itemize}
\item[i)] $\kappa=\sigma=0$ everywhere on $\mathcal M$
\end{itemize}
implies
\begin{itemize}
\item[ii)] $\Psi_0=\Psi_1=0$ everywhere on $\mathcal M$.  
\end{itemize}
\end{theorem} 

Third, using this theorem, and the definitions of $\Psi_0$ and
$\Psi_1$, we conclude that if our metric (\ref{cof}),
(\ref{choco}) satisfies (a) then, in our adapted coframe it has: 
\be
R_{2412}=R_{2424}=R_{2414}=R_{2423}=R_{2434}=0,\label{eiwei}
\ee
everywhere on $\mathcal M$. 

Fourth, we write down explicitly the second Cartan's structure
equations for indices $\{24\}$:
$$\der\Gamma_{24}+(-\Gamma_{12}+\Gamma_{34})\dz\Gamma_{24}=\tfrac12R_{24kl}\theta^k\dz\theta^l.$$
Due to (\ref{eiwei}) the r.h.s. of the above equality has only one
nonvanishing term:
$$\tfrac12R_{24kl}\theta^k\dz\theta^l=R_{2413}\theta^1\dz\theta^3.$$
Thus, assuming equations (a) we have an identity:
$$\der\Gamma_{24}+(-\Gamma_{12}+\Gamma_{34})\dz\Gamma_{24}-R_{2413}\theta^1\dz\theta^3=0.$$
Since $\Gamma_{24}$ as given by (\ref{wyco}) is a linear combination
of $\theta^1$ and $\theta^3$, after wedging this identity with
$\Gamma_{24}$, we conclude that:
\be\der\Gamma_{24}\dz\Gamma_{24}=0
\label{crf}
\ee  
everywhere on $\mathcal M$. A short calculation shows that
\be
\Gamma_{24}\dz\bar{\Gamma}_{24}=|\rho|^2\theta^1\dz\theta^2+\bar{\rho}\tau\theta^1\dz\theta^3-\rho\bar{\tau}\theta^2\dz\theta^3.\label{asub}
\ee
From now on we assume that $$\rho\neq 0$$ at every point on $\mathcal M$. This is equivalent to saying that the 
congruence of null and shearfree geodesics is diverging at points
where the associated CR structure has vanishing Levi form.

Our assumption about nonvanishing $\rho$, when compared with
(\ref{asub}), implies that 
\be
\Gamma_{24}\dz\bar{\Gamma}_{24}\neq 0\label{crf1}
\ee
at every point on $\mathcal M$.

Now we need the following lemma.  
\begin{lemma}\label{le}
Let $\varphi$ be a smooth complex valued 1-form defined locally in
$\bbR^n$, $n\geq 3$, such that $\varphi\dz\bar{\varphi}\neq 0$. Then
$$\der \varphi\dz\varphi \equiv 0\quad{\rm
if~ and~ only~ if}\quad \varphi=h\der\zeta$$
where $\zeta$ is a smooth complex function such that
$\der\zeta\dz\der\bar{\zeta}\neq 0$, and $h$ is a smooth nonvanishing complex
function.
\end{lemma}
\begin{proof}
Consider an open set
$U\in\bbR^n$ in which we have $\varphi$ such that
$\der\varphi\dz\varphi=0$ and $\varphi\dz\bar{\varphi}\neq 0$. We define {\it real} 1-forms
$\varphi^1={\rm Re}(\varphi)$ and $\varphi^2={\rm
  Im}(\varphi)$. They satisfy $\varphi^1\dz\varphi^2\neq 0$ in $U$. Our
assumption $\der\varphi\dz\varphi=0$, when written in terms of the real
1-forms $\varphi^1$, $\varphi^2$ is:
$$\der\varphi^1\dz\varphi^1-\der\varphi^2\dz\varphi^2+i(\der\varphi^2\dz\varphi^1+\der\varphi^1\dz\varphi^2)=0.$$
Taking the real and imaginary parts we have
\begin{eqnarray*}
&&\der\varphi^1\dz\varphi^1-\der\varphi^2\dz\varphi^2=0\\
&&\der\varphi^2\dz\varphi^1+\der\varphi^1\dz\varphi^2=0.
\end{eqnarray*}
Now the argument splits into two cases. For dimension $n\geq 4$ we
wedge 
the first equality, and then the second equality with
$\varphi^2$, and get
\be
\der\varphi^1\dz\varphi^1\dz\varphi^2=0,\quad\quad\der\varphi^2\dz\varphi^1\dz\varphi^2=0.\label{fre}
\ee
If $n=3$ these equations are trivially satisfied. 
In whatever dimension $n\geq 3$ we are in, once we have noticed that
(\ref{fre}) is true, we see that the \emph{real} forms 
$\varphi^1$, $\varphi^2$ form a closed differential ideal. Thus we can 
use the \emph{real} Fr\"obenius theorem, which  
implies that there exists a coordinate
chart $(x,y,u^\nu)$, $\nu=3,4,....n$, in $U$ such that 
$\varphi^1=t_{11}\der x+t_{12}\der y$
and $\varphi^2=t_{21}\der x+t_{22}\der y$, with some {\it real}
functions $t_{ij}$ in $U$ such that $t_{11}t_{22}-t_{12}t_{21}\neq
0$. Thus in the coordinates $(x,y,u^\nu)$ the form
$\varphi=\varphi^1+i\varphi^2$
can be written as $\varphi=c_1\der x+ c_2 \der
y$, where now $c_1$, $c_2$ are {\it complex} functions such that
$c_1\bar{c}_2-\bar{c}_1c_2\neq 0$ on $U$, so neither $c_1$ nor $c_2$
can be zero. The $\der\varphi\dz\varphi\equiv 0$
condition for $\varphi$ written in this representation is simply
$c_2^2\der(\tfrac{c_1}{c_2})\dz\der x\dz\der y\equiv 0$. Thus all the
partial derivatives $$\frac{\partial(\frac{c_1}{c_2})}{\partial u^\nu}\equiv 0,\quad\quad\forall\nu=3,4,...,n.$$ This means that
the ratio $\frac{c_1}{c_2}$ does not depend on $u^\nu$. This ratio defines
a nonvanishing {\it complex} function $F(x,y)=\frac{c_1}{c_2}$ of only
{\it two} real variables $x$ and $y$. Returning to
$\varphi$ we see that it is of the form $\varphi=c_2\big(\der y+F(x,y)\der
x\big)$. Consider the real bilinear symmetric
form $G=2\varphi\bar{\varphi}=|c_2|^2\Big(\der y^2+2\big(F(x,y)+\bar{F}(x,y)\big)\der
x\der y+|F(x,y)|^2\der
x^2\Big)$.
Invoking the
classical theorem on the existence of isothermal coordinates we are
able to find an open set $U'\subset
U$ with new coordinates $(\xi,\eta,u^\nu)$ in which $G=h^2(\der
\xi^2+\der \eta^2)$, where $h=h(\xi,\eta,u^\nu)$ is a real function in
$U'$. This means that in these coordinates $\varphi=h\der(
\xi+i\eta)=h\der\zeta$, and because of $\varphi\dz\bar{\varphi}\neq 0$
we have $\der\zeta\dz\der\bar{\zeta}\neq 0$. 
The proof in the other direction is obvious.
\end{proof}

Since the connection 1-form $\Gamma_{24}$ satisfies (\ref{crf}) and 
(\ref{crf1}), we can apply the above Lemma for
$\varphi=\Gamma_{24}$. Now, $n=4$ and we have
$$\Gamma_{24}=h\der\zeta\quad{\rm with}\quad\quad
\der\zeta\dz\der\bar{\zeta}\neq 0\quad{\rm on}\quad U'\subset\mathcal M.$$
Using (\ref{wyco}), which relates $\Gamma_{24}$ to
the coframe 1-forms $\theta^1$ and $\theta^3$, and expressing these
two in terms of the 1-forms 
$(\lambda,\mu,\bar{\mu})$ by (\ref{choco}), we get:
$$h\der\zeta=\Gamma_{24}=-P(\bar{\rho}\mu+\bar{\tau}\lambda).$$
Since the function $P$ is nowhere vanishing we can write this last
equation as
$-P^{-1}h\der\zeta=\bar{\rho}\mu+\bar{\tau}\lambda$. Wedging this
with $\der\zeta\dz\lambda$ we get
$\bar{\rho}\der\zeta \dz\lambda\dz\mu=0$ on $U'$. 
Because of our assumption that $\rho$ is nowhere vanishing we finally obtain \be\der\zeta
\dz\lambda\dz\mu=0\quad{\rm with}\quad\der\zeta\dz\der\bar{\zeta}\neq
0\quad{\rm on}\quad U'\subset\mathcal M.\label{koty}\ee
The last equation pullsback to the CR manifold $M$ providing a CR
function there.

Our construction of $\zeta$ obviously works if the CR structure is of
class $C^3$. Actually we think that
class $C^{2,1}$ suffices, see \cite{ht}.\\

The last part of the proof consists in giving a geometric
interpretation to the Einstein conditions (a).\\

To discuss this we go back to a general spacetime
$\mathcal M$ equipped with a null congruence associated with a
vector field $k$. We do not require that this congruence is geodesic
and shearfree here. Such a congruence defines a class of adapted
coframes (\ref{cof}), (\ref{cofr}). It follows from equations
(\ref{cofr}) that although $k$ does not specify 
the coframe 1-forms $\theta^1$, $\theta^3$ uniquely, we have
$${\theta^1}'\dz{\theta^3}'=A{\rm e}^{i\phi}\theta^1\dz\theta^3.$$
This may be used to define complex valued vector fields $Y$ on 
$\mathcal M$ such that 
$$Y\hook (\theta^1\dz\theta^3)=0.$$
The complex valued distribution $N$ consisting of all such vector
fields is uniquely defined by the null congruence on $\mathcal M$. It
follows that
$$N=\{a e_2+b e_4\}$$
where $a,b$ are arbitrary complex valued functions on $\mathcal M$ and
$e_2$, $e_4$ is a part of the null frame $(e_1,e_2,e_3,e_4)$ dual
to $(\theta^1,\theta^2,\theta^3,\theta^4)$,
i.e. $e_i\hook\theta^j=\delta^j_{~i}$ and, in particular, $e_4\dz k=0$. 

Thus a null congruence defines at each point $x\in M$ a 
2-complex-dimensional plane $N_x$. These planes are called
$\alpha$ planes. They have the property of being \emph{totally null}:
$$g(Y_1,Y_2)=0,\quad\quad\forall Y_1,Y_2\in N.$$
Thus all vectors in $N_x$ are null and orthogonal to each other.

The Ricci tensor of a spacetime may be
considered as a symmetric bilinear form. We extend it 
to the complexification by linearity. 

Now combining these two facts we may require that we have a spacetime
in which the  Ricci tensor, extended to the complexification, has a
similar property with respect to $N$ as the metric. We say that
the complexified Ricci tensor $Ric(g)$ \emph{vanishes on the  
distribution} $N$ iff 
$$Ric(g)(Y_1,Y_2)=0,\quad\quad\forall Y_1,Y_2\in N.$$
We will denote this condition by $Ric(g)_{|N}=0$. 
Obviously it is weaker than the Ricci flatness condition. Actually if
$N$ is the distribution of $\alpha$ planes associated with
the congruence, then in an adapted coframe
$(\theta^1,\theta^2,\theta^3,\theta^4)$
$$Ric(g)_{|N}=0\quad\quad\Leftrightarrow\quad\quad R_{22}=R_{24}=R_{44}=0.$$
Hence the vanishing of Ricci on $N$ is precisely equivalent to our conditions (a).

This completes the proof of Theorem \ref{firstcr}.

\subsection{Reduction of the Einstein equations to the CR manifold}\label{ree}
In this section we derive a maximally reduced system corresponding to the Einstein
equations for a spacetime admitting a congruence of
null and shearfree geodesics. We will assume
that the congruence has nonvanishing twist at every point of the spacetime. This is the
same as to assume that the underlying CR structure is strictly
pseudoconvex.

As we know, the Einstein equations for such spacetimes split
conveniently into three types of equations, which in Section
\ref{sedf} were denoted by (a), (b) and (c). Our reduction procedure will
follow this split: we first impose equations (a), then (b) and finally
(c).

We start with equations (a). They will be reduced according to the
following scheme. As we have proven in Section \ref{sedf} equations (a)
imply that the corresponding CR structure admits at least one
CR-function. This result enables us to start our reduction procedure
with a 4-manifold ${\mathcal M}=\bbR\times M$ by taking the metric in the form
\be
g=2(\theta^1\theta^2+\theta^3\theta^4),\label{mde}
\ee
where 
\begin{eqnarray*}
&&\theta^1=P\mu,\quad\quad\theta^2=P\bar{\mu},\quad\quad\theta^3=P\lambda,\quad\quad\theta^4=P[\der r+W\mu+\bar{W}\bar{\mu}+H\lambda],
\end{eqnarray*}
with the 1-forms $(\lambda,\mu,\bar{\mu})$ satisfying
\begin{eqnarray}
&&\der\mu=0,\quad\quad\quad\der\bar{\mu}=0,\label{onecr}\\
\der\lambda&=&i\mu\dz\bar{\mu}+(c\mu+\bar{c}\bar{\mu})\dz\lambda.\label{psec}
\end{eqnarray}
The forms $(\lambda,\mu,\bar{\mu})$ define the CR structure on $M$. Note 
that formula (\ref{onecr}) follows from our result on the existence of one CR
function. This result enables us to put  
$$\mu=\der\zeta,$$
where $\zeta$ is the first CR function obtained in (\ref{koty}). We
remain with this choice for $\mu$ in the rest of this section. Formula (\ref{psec}), which says that there is a choice of
$\lambda$ such that the coefficient of the $i\mu\dz\bar{\mu}$ term is
equal to one, is equivalent to our assumption about strict 
pseudoconvexity. 

At this stage we introduce a basis
$$(\partial_0,\partial,\bar{\partial})$$ of vector fields on $M$, which
is dual to the coframe $(\lambda,\mu,\bar{\mu})$ on
$M$. The complex vector field $\bar{\partial}$ is the \emph{tangential
CR operator on $M$}. 

We also note that the closure of the system (\ref{onecr})-(\ref{psec})
implies that:
\be
\partial\bar{c}=\bar{\partial}c,\label{ee0}
\ee   
so that $\partial\bar{c}$ is \emph{real}.

Thus from now on we assume that we have a CR manifold
$M$ with the defining forms $(\lambda,\mu,\bar{\mu})$ satisfying 
(\ref{onecr})-(\ref{psec}). Our goal is to lift this CR structure to
an Einstein spacetime. 
\begin{remark}
We stress that although we did not impose the full equations (a) on
our spacetime, we already used in (\ref{onecr})-(\ref{psec}) a
consequence 
of these equations, which
enabled us to assume that $(\lambda,\mu,\bar{\mu})$ are in the form
(\ref{psec}). It is a justified procedure: since ultimately we are
interested in the maximal reduction of the full system (a), we may freely use
its consequence at any stage of the reduction procedure.
\end{remark}
The reduction of equations (a) goes as follows:
\begin{itemize}
\item We first use the Goldberg-Sachs theorem, which says that if
  equations (a) are satisfied then $R_{2412}+R_{2434}=0$. Modulo 
  complex conjugation this is equivalent to the requirement that the
  function $W$ satisfies:
$$W_r-iW_{rr}=0.$$
This equation may be easily integrated, proving that the most general
form of a $W$ satisfying (a) is given by
\be
W=i {\rm e}^{-ir}x+y,\label{ee1}
\ee
where the complex functions $x$ and $y$ are
$r$-independent, $x_r=y_r=0$. Thus the requirement that our
spacetime is algebraically special (the requirement that is implied by
(a)) is equivalent to the form (\ref{ee1}) of the function $W$.
\item The first of equations (a), namely $R_{44}=0$, is equivalent to the
  differential equation on $P$:
$$-4PP_{rr}+P_r^2+P^2=0.$$
This again can be easily solved to get:
\be
P=\frac{p}{\cos(\frac{r+s}{2})},\label{ee2}
\ee
with real functions $p\neq 0$ and $s$ satisfying $p_r=s_r=0$.
\item Equation $R_{24}=0$ is equivalent to 
$$
y=-i c-2i\partial\log p+\partial s+2i{\rm e}^{is}x,
$$   
i.e. to an equation which expresses the function $y$ in terms of the
unknowns $p,s,x$ and the CR quantities $c$ and $\partial$. 
\item Now it is convenient to introduce a \emph{new} unknown $t$, a 
  complex valued function on $M$, which replaces the unknown $x$. The
  quantity $x$ is related to $t$ via:
\be
{\rm e}^{is}x=c+2\partial\log p -t.\label{ee3}
\ee
This enables us to write $y$ as
\be
y=ic+2i\partial\log p+\partial s-2i t.\label{ee4}
\ee
\item In terms of $t$ the Einstein equation $R_{22}=0$ is
  equivalent to 
\be
\partial t+(c-t)t=0.\label{ee5}
\ee
\end{itemize} 
We sumarize this in the following proposition. 
\begin{proposition}\label{proa}
The metric (\ref{mde})-(\ref{psec}) satisfies the Einstein equations
$R_{44}=R_{24}=R_{22}=0$ if and only if the function $W$ is given by
(\ref{ee1}), the function $P$ is given by (\ref{ee2}) with 
$p, s$ (real), $x,y,t$ (complex) satisfying $p_r=s_r=x_r=y_r=t_r=0$, 
definitions (\ref{ee3}), (\ref{ee4}) and the differential equation
(\ref{ee5}). 
\end{proposition} 
Observe that equation (\ref{ee5}) always has the solution $t\equiv
0$. Using this observation we prove the converse of 
Theorem \ref{firstcr} in the strictly pseudoconvex case. Indeed if we
have a strictly pseudoconvex 3-dimensional CR manifold with one  
CR function $\zeta$ such that $\der\zeta\dz\der\bar{\zeta}\neq
0$, we define $\mu=\der\zeta$. We next choose $\lambda$ so that (\ref{psec})
is satisfied. Then we take $t\equiv 0$ and choose sufficiently smooth
arbitrary real functions $p$ and $s$. Using them we define a lift to a 
spacetime having the metric $g$ as in (\ref{mde}), with $P$
given by (\ref{ee2}), and $W$ given by (\ref{ee1}), (\ref{ee3}), (\ref{ee4}). Then it follows that such a metric satisfies Einstein's
equations (a), independently of the choice of a real function $H$ on
$\mathcal M$. Combining this with Theorem \ref{firstcr} we proved
Theorem \ref{fistcr}. \\

We now pass to the reduction of equations (b). The steps here are:
\begin{itemize}
\item Assuming equations (a) to be satisfied, so that the metric is
  given by Proposition \ref{proa}, we first impose a consequence of
  (b), namely $R_{12}+R_{34}=2\Lambda={\rm const}$, which is just the
  condition that the Ricci scalar is constant and equal to $4\Lambda$. This equation
  determines the $r$ dependence of the function $H$. It reads:
\be
H=\frac{m}{p^4}{\rm e}^{2i(r+s)}+\frac{\bar{m}}{p^4}{\rm
  e}^{-2i(r+s)}+Q{\rm e}^{i(r+s)}+\bar{Q}{\rm e}^{-i(r+s)}+h,\label{ee6}
\ee 
with
\begin{eqnarray*}
Q&=&\frac{3m+\bar{m}}{p^4}+\frac23\Lambda p^2+\frac{2\partial p\bar{\partial}p-p(\partial\bar{\partial}p+\bar{\partial}\partial
  p)}{2p^2}-\\
&&\frac{i}{2}\partial_0\log
  p-2t\bar{\partial}p-\bar{t}\partial\log
  p+\frac32\bar{\partial}t-\bar{c}t-\frac12c\bar{t}+\frac52
  t\bar{t}+\partial\bar{t}-\bar{\partial}c\\&&\\
h&=&3\frac{m+\bar{m}}{p^4}+2\Lambda p^2+\frac{2\partial
  p\bar{\partial}p-p(\partial\bar{\partial}p+\bar{\partial}\partial
  p)}{p^2}-\\&&3(t\bar{\partial}\log p+\bar{t}\partial\log
  p)+\frac52(\partial\bar{t}+\bar{\partial}t)+6t\bar{t}-\frac32(c\bar{t}+\bar{c}t)-2\bar{\partial}c+\partial_0 s.
\end{eqnarray*}
Note that $h$ is \emph{real} due to (\ref{ee0}). The unknown $m$
is complex and satisfies $m_r=0$.  
\item At this stage we have found the explicit $r$ dependence of the entire
  metric (\ref{mde})-(\ref{psec}). For the full determination of this
  dependence we only needed equations (a) and the subset of equations (b)
  given by $R_{12}+R_{34}=2\Lambda={\rm const}$. 
\item We now impose another consequence of (b), namely
  $R_{12}-R_{34}=0$. This, together with $R_{12}+R_{34}=2\Lambda={\rm
  const}$, is equivalent to (b). The reduction of this equation gives
  the following differential equation connecting $p$, $t$ and $m$:
\begin{eqnarray}
&&[~\partial\bar{\partial}+\bar{\partial}\partial
   +\bar{c}\partial+c\bar{\partial}+\tfrac12c\bar{c}+\tfrac34(\partial\bar{c}+\bar{\partial}
    c)-\tfrac32(\partial\bar{t}+\bar{\partial}t+t\bar{t})~]p=\label{ee7}\\
&&\frac{m+\bar{m}}{p^3}+\tfrac23\Lambda p^3.\nonumber
\end{eqnarray}    
This completes the reduction of equations (b). 
\end{itemize}
We summarize in the following theorem. 
\begin{theorem}
A strictly pseudoconvex CR structure $(M,(\lambda,\mu))$ 
lifts to a spacetime satisfying
Einstein equations (a) and (b) if and only if it admits at least one CR
function $\zeta$ with $\der\zeta\dz\der\bar{\zeta}\neq 0$ and, in
addition, it admits a solution to equation (\ref{ee7}) for a
real function $p$ on $M$, with $c,t$ obeying respectively 
(\ref{psec}) and (\ref{ee5}). 
\end{theorem}
If we do not insist on the full system (a) and (b) we conclude the
following remarkable theorem.  
\begin{theorem}\label{psyk}
A strictly pseudoconvex CR structure $(M,(\lambda,\mu))$ 
lifts to a spacetime having constant Ricci scalar and satisfying
equations (a) if and only if it admits at least one CR
function $\zeta$ with
$\der\zeta\dz\der\bar{\zeta}\neq 0$. In such case the spacetime metric
satisfying (a) and having constant Ricci scalar equal to $4\Lambda$ is 
(\ref{mde})-(\ref{psec}) with $W$ given by (\ref{ee1}), $P$ given by
(\ref{ee2}), $x,y$ given by (\ref{ee3}), (\ref{ee4}), and $H,Q,h$ given by 
(\ref{ee6}). The functions $m$ (complex), $p$ (real) are arbitrary, and
the complex function $t$ satisfies the partial differential equation
(\ref{ee5}). 
\end{theorem}
This theorem is remarkable for the reasons highlighted in the
following remarks.

\begin{remark}
Given a CR structure with one CR
function, to determine the most general lift to a spacetime with a metric satisfying the Ricci
conditions of Theorem (\ref{psyk}), we need to have a general
solution for only \emph{one} complex equation
  (\ref{ee5}) for the 
  complex function $t$ on the CR manifold. This equation has always
  one solution, namely $t=0$. Surprisingly the question if this
  equation 
has other solutions is equivalent to the question if the CR structure
  admits more CR functions, and hence is locally embeddable. 
To see this take 
the CR manifold $(M,(\lambda,\mu))$ satisfying
  (\ref{onecr})-(\ref{psec}) and consider the complex 1-form $\varphi$ 
   defined by:
\be
\varphi=\mu+i \bar{t}\lambda.\label{0psy}\ee
Then due to (\ref{onecr})-(\ref{psec}) we have 
$$\der\varphi\dz\varphi=i~[~\bar{\partial}\bar{t}+(\bar{c}-\bar{t})\bar{t}~]~
\mu\dz\bar{\mu}\dz\lambda.$$
Thus $\der\varphi\dz\varphi=0$ is equivalent to the Einstein equation
(\ref{ee5}). Since obviously $\varphi\dz\bar{\varphi}\neq 0$, then according
to Lemma \ref{le}, $\phi$ defines a CR function $\eta$ such that
  $h\der\eta=\varphi$. Thus the equation \be
h\der\eta=\mu+i\bar{t}\lambda\label{psy}\ee
relates the CR functions $\eta$
and solutions $t$ of the Einstein equation 
  (\ref{ee5}). As an example take the trivial solution $t=0$ of
(\ref{ee5}). Since $\mu=\der\zeta$ and since for $t=0$ equation
(\ref{psy}) gives $h\der\eta=\mu$,  the CR
function $\eta$ is dependent on the CR function $\zeta$
of (\ref{koty}). To get a $\zeta$-independent CR function $\eta$ we
need a nonzero solution of (\ref{ee5}). That this requirement is also
sufficient follows from the relation
$h\der\eta\dz\der\zeta=i\bar{t}\lambda\dz\mu$ implied by (\ref{psy}). This proves the following corollary.
\begin{corollary}\label{2psy}
Every nonzero solution $t$ of the Einstein equation (\ref{ee5})
provides a second CR function $\eta$ such that
$\der\eta\dz\der\zeta\neq 0$. Also
the converse is true: every CR function $\eta$ 
defines a complex function $t$ satisfying the Einstein equation
(\ref{ee5}). The transformation between $\eta$ and $t$ is given by 
$h\der\eta=\der\zeta+i\bar{t}\lambda$.   
\end{corollary} 
Thus, in particular, if the
  starting CR structure is locally embedded, which means that the
  general solution to the tangential CR equation is explicitly known,
  we may write the general solution for $t$ satisfying (\ref{ee5}) and
  obtain the most general solution for the metric satisfying 
the highly nonlinear system of
  equations (a) and $R_{12}+R_{34}=2\Lambda={\rm const}$.  
\end{remark}
\begin{remark}
Another remarkable feature of Theorem \ref{psyk} is that the metrics
(\ref{mde})-(\ref{psec}) satisfying (a) and $R_{12}+R_{34}=2\Lambda={\rm
  const}$ have an explicit $r$ dependence which is very
particular. Note that the functions $W$ of (\ref{ee1}) and $H$ of
(\ref{ee6}) are \emph{periodic} in $r$ with period $2\pi$. If we
forget about the conformal factor we may write the metrics in the form
\cite{nurphd} 
\begin{eqnarray}
\hat{g}&=&2p~\Big[~\mu\bar{\mu}~+~\lambda~[~\der r ~+~\big(i {\rm e}^{-ir}x+y\big)\mu~+~\big(-i {\rm e}^{ir}\bar{x}~+~\bar{y}\big)\bar{\mu}~+~\nonumber\\&&\big(~\frac{m}{p^4}{\rm e}^{2i(r+s)}+\frac{\bar{m}}{p^4}{\rm
  e}^{-2i(r+s)}+Q{\rm e}^{i(r+s)}+\bar{Q}{\rm
  e}^{-i(r+s)}+h~\big)\lambda]~\Big],\label{fefer}
\end{eqnarray}
which is periodic and \emph{regular} in $r$. Thus, what the Einstein
equations (a) and $R_{12}+R_{34}=2\Lambda={\rm const}$ impose on the
spacetime ${\mathcal M}=\bbR\times M$ is a \emph{circle bundle
  structure} $\bbS^1\to\hat{\mathcal M}\to M$ 
on the Lorentzian manifold $(\hat{\mathcal M},\hat{g})$
which has $\mathcal M$ as its universal cover. The Lorentzian manifold 
$(\hat{\mathcal M},\hat{g})$ is called by physicists 
a \emph{conformal (fiberwise) compactification} of $({\mathcal M},g)$. It is used by them
to study the asymptotic behavior of a gravitational field. We
summarize in the following corollary.
\begin{corollary}
The Einstein equations (a) and the constancy of the Ricci scalar imposed on
the metrics (\ref{mde})-(\ref{psec}) imply that all the metric
functions are periodic in the $r$ coordinate, so that there is a natural circle
bundle over the strictly pseudoconvex CR manifold onto which all the
Einstein metrics (\ref{mde})-(\ref{psec}) descend. 
\end{corollary}  
\end{remark}
\begin{remark}
In 1976 Fefferman \cite{fefe} introduced a natural \emph{conformal} Lorentzian
metric $\hat{g}_F$ on a circle bundle $\bbS^1\to \hat{\mathcal M}_F\to
M$ over any strictly pseudoconvex 3-dimensional
CR manifold $(M,(\lambda,\mu))$ embedded in $\bbC^2$. A natural question is how
his circle bundle and his
Lorentzian metrics are related to our $(\hat{\mathcal M},\hat{g})$
above. The answer is the following:
\begin{itemize}
\item Fefferman metrics constitute a simple subclass of our metrics 
(\ref{fefer}).
\item They happen to be conformally Einstein \emph{only} in the case when the corresponding CR manifold is locally CR equivalent to the Heisenberg group CR structure  \cite{lew}; in such case $\hat{g}_F$ is conformally flat.
\item Given a CR structure as in (\ref{onecr})-(\ref{psec}) the
  Fefferman metric $\hat{g}_F$ is obtained from our $\hat{g}$ by
  putting $x=m=Q=0$, $y=-\frac{i}{3}c$,
  $h=-\frac{1}{12}(\partial\bar{c}+\bar{\partial}c)$.
\item Thus, in our setting, the Fefferman metrics (or, strictly speaking, their generalizations to strictly pseudoconvex CR manifolds which admit one CR function) are represented by
\be
\hat{g}_F=2~\Big[~\mu\bar{\mu}~+~\lambda~[~\der r
    ~-~\tfrac{i}{3}c\mu~+~\tfrac{i}{3}\bar{c}\bar{\mu}~-~\tfrac{1}{12}(\partial\bar{c}+\bar{\partial}c)~\big)\lambda]~\Big].\label{ffe}
\ee
\end{itemize}
Note that the Fefferman metrics are $r$-independent. This reflects the
well known fact \cite{lew,nurphd,nurple,spar} that the null congruence of shearfree geodesics
associated with the $k=\partial_r$ direction is a \emph{conformal
  Killing} vector for each Fefferman metric. Actually, the above
formula for the Fefferman metric for CR manifolds having one CR
function is obtained by (i) imposing the requirement that the metric (\ref{fefer})
has a conformal Killing vector alligned with $k=\partial_r$ (this
forces $x$, $m$ and $Q$ to vanish), and (ii) imposing another
requirement that the metrics (\ref{fefer}) are of Petrov type $N$ (this
specifies that $y$ and $h$ must be expressed in terms of $c$ as
above). That the requirements $(i)$ and $(ii)$ are neccessary and
sufficient to distinguish the Fefferman metrics among metrics
(\ref{fefer}) is a well known fact \cite{lew,nurphd,nurple,spar}.

Looking at the Fefferman metrics (\ref{ffe}) one may say that the
circle bundle structure of the spacetime is not visible, since the
metrics are \emph{constant} along the $r$-direction. 
Only if a more general class of metrics (\ref{mde})-(\ref{psec}) 
is taken into account does the
circle bundle structure emerge. And it emerges in a natural way, 
as a consequence of the Einstein
equations (a) and $R_{12}+R_{34}=2\Lambda={\rm const}$.

The fact that the Fefferman metrics are of Petrov type $N$ everywhere 
essentially means that the Weyl tensor has only one nonvanishing
complex component. This is proportional to 
$$\Psi_4=\partial\bar{\partial}\partial c+3c\bar{\partial}\partial
c-7ic\partial_0c-3i\partial\partial_0c+(\partial
c+2c^2)\bar{\partial}c.$$
The vanishing or not of $\Psi_4$ is a CR invariant property. It is
also a conformal invariant property. In the CR context, vanishing of
$\Psi_4$ is an if and only if condition for the CR structure
(\ref{onecr})-(\ref{psec}) to be CR equivalent to the Heisenberg group
CR structure \cite{BDS,lew}. In the conformal context, the vanishing of $\Psi_4$ means
that $\hat{g}$ is conformally flat. Thus the $\Psi_4$ component of the
Weyl tensor for $\hat{g}_F$ is proportional to Cartan's lowest order
invariant \cite{carcr} of a strictly pseudoconvex CR structure.

We thus have the following corollary. 
\begin{corollary}
The CR structure (\ref{onecr})-(\ref{psec}) is locally CR equivalent to
the Heisenberg group CR structure if and only if the function $c$ in
(\ref{psec}) satisfies  
$$\partial\bar{\partial}\partial c+3c\bar{\partial}\partial
c-7ic\partial_0c-3i\partial\partial_0c+(\partial
c+2c^2)\bar{\partial}c=0.$$
\end{corollary}  
\end{remark}

\begin{remark}
We are now prepared to give the geometric interpretation of the
leading terms in the partial differential equation (\ref{ee7}). Having
chosen a strictly pseudoconvex CR structure
$(M,(\lambda,\mu))$ and functions $t$ and $m$ on $M$ this
is a partial differential equation for a real function $p$ on $M$. Below we give
the interpretation of the \emph{linear} operator on its left hand side.

Since the Fefferman metrics $\hat{g}_F$ are defined up to a conformal scale 
it is reasonable to consider the \emph{conformally invariant} 
wave operator $(-*\der *\der-\tfrac16 R)$, 
with $-*\der *\der$ being the D'Alembert operator and $R$ being the Ricci
scalar; both being calculated in the Fefferman metric $\hat{g}_F$. Let us apply this operator
to a real function $f$ on the Fefferman bundle $\hat{\mathcal M}_F$,
which is constant along the fibres. Remarkably we get \cite{nurphd}:
\be
(-*\der *\der-\tfrac16
R)f~=~[~\partial\bar{\partial}+\bar{\partial}\partial+c\bar{\partial}+\bar{c}\partial+\tfrac12
  c\bar{c}+\tfrac38(\partial\bar{c}+\bar{\partial}c)~]f.\label{kg}
\ee
It is worthwhile to add that 
$R=-3
c\bar{c}-\tfrac94(\partial\bar{c}+\bar{\partial}c)$ and that $*\der
*\der
f=-(\partial\bar{\partial}+\bar{\partial}\partial+c\bar{\partial}+\bar{c}\partial)f$, 
but neither $R$ nor the 3-dimensional operator
$\partial\bar{\partial}+\bar{\partial}\partial+c\bar{\partial}+\bar{c}\partial$
has a CR geometrical meaning. Only their sum
$$\triangle_{CR}=\partial\bar{\partial}+\bar{\partial}\partial+c\bar{\partial}+\bar{c}\partial+\tfrac12
c\bar{c}+\tfrac38(\partial\bar{c}+\bar{\partial}c)$$
is CR meaningful.

Using (\ref{kg}) we rewrite the Einstein equation (\ref{ee7}) in 
the following equivalent form:
$$\Big[\triangle_{CR}+\tfrac38(\partial\bar{c}+\bar{\partial}
    c)-\tfrac32(\partial\bar{t}+\bar{\partial}t+t\bar{t})\Big]p=\frac{m+\bar{m}}{p^3}+\tfrac23\Lambda p^3.$$
The appearence of the $\tfrac38(\partial\bar{c}+\bar{\partial}
    c)$ term here in the potential is unpleasent, but we can not avoid it.
\end{remark}

Finally we pass to equations (c):
\begin{itemize}
\item Assuming (a) and (b) are satisfied, we first reduce the complex equation
  $R_{13}=0$. This is equivalent to one complex equation 
\be
\partial m+3(c-t)m=0,\label{ee8}
\ee
for the complex function $m$ of (\ref{ee6}).
\item Now assuming that the metric (\ref{mde})-(\ref{psec}) satisfies
  the Einstein equations (\ref{ee1}), (\ref{ee2}), (\ref{ee3}),
  (\ref{ee4}), (\ref{ee5}), (\ref{ee6}), (\ref{ee7}) and (\ref{ee8})
  we calculate the Weyl tensor. Since, via Goldberg-Sachs, the metric
  is algebraically special, only three of the five Weyl scalars $\Psi_0$,
  $\Psi_1$, $\Psi_2$, $\Psi_3$, $\Psi_4$ are in principle 
  nonvanishing. These are: $\Psi_2$, $\Psi_3$, $\Psi_4$. The Weyl
  spinor $\Psi_2$, whose vanishing means that the metric is of Petrov
  type $III$, $N$ or $0$, has a particularly neat form:
\be
\Psi_2=\frac{(1+{\rm e}^{i(r+s)})^3}{2p^6}m.\label{tyd}
\ee     
\item Although we succeded in reducing the last equation $R_{33}=0$,
  the explicit reduced form of it is too complicated to be presented here. 
\end{itemize}
The following remark is in order. 
\begin{remark}\label{rty}
If a lift to the spacetime $({\mathcal M},g)$ 
of \emph{any} CR structure $(M,(\lambda,\mu))$ 
is considered, one
can try to write down curvature conditions that are compatible with
the underlying CR geometry of the associated congruence of shearfree
and null geodesics. We already met the curvature conditions that
respect the underlying CR geometry. These are the Einstein equations
(a), or in more geometric terms, conditions forcing the complexified
Ricci tensor to be identically zero on the associated distribution of
$\alpha$ planes. It turns out that the Einstein conditions (a), (b) and
$R_{13}=0$ also have geometric meaning. They are equivalent to
\be
Ric(g)=\Lambda g+\Phi\lambda\otimes\lambda,\label{eep}
\ee
where $\Phi$ is an arbitrary real function on $\mathcal M$. Physicists
call these equations the \emph{cosmological constant Einstein's
  equations with a pure radiation field}, since they describe 
gravitational fields with the energy momentum tensor in which all the
energy is propagated with the speed of light along the direction
determined by the congruence of shearfree geodesics defined by $M$.  
\end{remark}
\subsection{The second CR function}\label{reee}
We start this section by considering a 3-dimensio\-nal CR manifold which
lifts to a spacetime with a twisting
congruence of null and shearfree geodesics and which satisfies
Einstein equations (a).  
Theorem \ref{firstcr} assures that such a CR manifold has at least one
CR function, say $\zeta$, with 
$\der\zeta\dz\der\bar{\zeta}\neq 0$. Our approach to the problem of
obtaining an independent CR function, say $\eta$, is by finding a complex equation, let us call it
($2^{{\rm nd}}CR$),
which when assumed to be satisfied, guarantees that $\eta$ exists. 
The main idea here is to find 
($2^{{\rm nd}}CR$) among the full system
of Einstein equations (a), (b) and (c), especially among the reduced
Einstein's equations of Section \ref{ree}.

It turns out that, depending on some additional assumptions about the
lifted spacetime, various choices of $(2^{{\rm nd}}CR$) are possible. Among all
of these choices the simplest is to consider equation (\ref{ee5}) as
the ($2^{{\rm nd}}CR$). Indeed, if our CR manifold lifts to a
spacetime whose metric (\ref{mde})-(\ref{psec}) has nonvanishing $t$ 
in (\ref{ee5}) then, as we already noticed in Corollary (\ref{2psy}),
it is locally embeddable, with an embedding given by means of the CR functions 
$\zeta$ of (\ref{koty}) and $\eta$ of (\ref{psy}). The trouble with
equation (\ref{ee5}) is that the Einstein equations (a) do not
guarantee that (\ref{ee5}) has any other solution than $t\equiv
0$. Nevertheless equation (\ref{ee5}) may be used as a \emph{criterion for
embeddability}. Suppose one knows that
\begin{itemize}
\item[i)] a 3-dimensional strictly pseudoconvex CR manifold admits a
  CR function $\zeta$ such that $\der\zeta\dz\der\bar{\zeta}\neq 0$
  and, in addition, one knows that
\item[ii)] this manifold lifts to a spacetime with a metric $g$ 
which satisfies Einstein's equations (a).
\end{itemize} 
Then he can write the metric $g$ in the form (\ref{mde})-(\ref{psec})
with $\mu=\der\zeta$ and simply \emph{calculate} the function $t$. If he
finds that $t\neq 0$, then he concludes that the CR structure is
locally embeddable. This is due to the fact that the calculated $t$
automatically satisfies (\ref{ee5}) since $g$ satisfies (a). 

\begin{remark}
At this stage we remark that equation (\ref{ee5}) is interesting on
its own, without any reference to the fact that it originates from the
Einstein equation for the lifted spacetime. Indeed, it follows from
our discusion above, that one can use this equation to get a sharp
criterion for the embeddability of a strictly pseudoconvex CR structure
that admits one CR function. Here the procedure is as follows:
\begin{itemize}
\item Suppose we are given a strictly pseudoconvex 3-dimensional
  CR manifold $M$ which has one CR function $\zeta$ such that
  $\der\zeta\dz\der\bar{\zeta}\neq 0$.
\item Given $\zeta$, we write $\mu=\der\zeta$, and choose $\lambda$ so that
  equation (\ref{psec}) is satisfied. This, in particular, defines the
  function $c$ on $M$. 
\item Define $(\partial,\bar{\partial},\partial_0)$ as dual to
  $(\mu,\bar{\mu},\lambda)$.
\item Consider the equation $\partial t+(c-t)t=0$ 
for a complex function $t$ on $M$.
\item Then we have the following theorem.  
\end{itemize} 
\begin{theorem}
The CR structure $(M,(\lambda,\mu=\der\zeta))$ is locally embeddable 
if and only
if the equation $\partial t+(c-t)t=0$ has a solution $t$ such that $t\neq 0$.
\end{theorem}
We may combine this result with a result of Hanges \cite{hong}, who
found another criterion for the existence of the second CR function for a
3-dimensional strictly pseudoconvex CR manifold. It is well known,
that if a 3-dimensional strictly pseudoconvex 
CR manifold $M$ admits one CR function $\zeta$ as above, then one 
can supplement $\zeta$ and $\bar{\zeta}$ by a real function $u$ on $M$ so that
$({\rm Re}\zeta,{\rm Im}\zeta,u)$ constitute a coordinate system on
$M$, in which $$\mu=\der\zeta,\quad\quad {\rm
  and}\quad\quad\lambda=\frac{\der
  u+L\der\zeta+\bar{L}\der\bar{\zeta}}{i(\bar{\partial}L-\partial\bar{L})},$$
and in which the complex valued function $L=L(\zeta,\bar{\zeta},u)$
vanishes at the origin, $L(0,0,0)=0$. Hanges'
result is that the CR structure $(M,(\lambda,\mu))$ is locally
embeddable near the origin if
and only if the function $L=L(\zeta,\bar{\zeta},u)$ 
is the \emph{boundary value} of a function
$\tilde{L}=\tilde{L}(\zeta,\bar{\zeta},w)$ which is \emph{holomorphic} 
in the complex variable $w=u+iv$. Using this
result and writing the differential equation $\partial t+(c-t)t=0$ in
the local coordinates $(\zeta,\bar{\zeta},u)$ we get the following remarkable corollary. 
\begin{corollary}
The nonlinear partial differential equation 
$$(\partial\bar{L}-\bar{\partial}L)\partial
t-[\partial(\partial\bar{L}-\bar{\partial}L)+(\partial\bar{L}-\bar{\partial}L)(L_u+t)]t=0,$$
with $\partial=\partial_{\zeta}-L\partial_u$, and with
the complex valued function $L=L(\zeta,\bar{\zeta},u)$ vanishing at the
origin, is locally solvable near the origin for
a complex function $t\neq 0$ if and only if $L$ 
is the boundary value of a function
$\tilde{L}=\tilde{L}(\zeta,\bar{\zeta},w)$ which is holomorphic in the
variable $w$. If this is the
case the CR structure $(M,(\lambda,\mu))$ is locally embeddable.
\end{corollary} 
\end{remark}

Returning to our discussion of the relations between the second CR
function and the Einstein equations of the lifted spacetime, we are
now in a position to
say that, if we have a solution $t\neq 0$ of equation (\ref{ee5}), the
problem of the local embeddability of a 3-dimensional manifold $M$ is
solved. If we are in an unlucky situation which negates the existence
of a solution $t\neq 0$, two things may happen:
\begin{itemize}
\item either the only solution to (\ref{ee5}) is $t$ vanishing everywhere,
\item or $t=0$ at a point around which we want to embedd $M$ into $\bbC^2$.
\end{itemize}

In the first case, we can put $t\equiv 0$ in all the equations we have 
derived in Section \ref{ree}. In the second case, some care is needed,
and we need some preparations. In what follows, our considerations will be of a bit 
more general nature than is required to treat this case, but at a
certain moment, they will lead us to the conclusion that the case $t=0$ at a
point is, actually, the same as the case $t\equiv 0$. 

To get to this conclusion, consider a situation in which we have a
3-dimensional CR structure $(M,(\lambda,\mu))$ with 
$\mu$ and $\lambda$ satisfying (\ref{onecr})-(\ref{psec}). Assume, in
addition, that the CR structure admits complex functions $h\neq 0$ and
$t_0$ such that the complex valued 1-form \be
\mu'=h^{-1}(\mu+i\bar{t}_0\lambda)\quad\quad {\rm is~closed},\quad\quad
\der\mu'=0.\label{323a}\ee 
We \emph{assume} it is the case. If we have such $h$ and
$t_0$, we define 
\be\lambda'=|h|^{-2}\lambda.\label{lap}\ee 
The forms 
$(\lambda',\mu',\bar{\mu}')$ define the same CR structure on $M$ as
the forms $(\lambda,\mu,\bar{\mu})$. Moreover, because of our 
choice of $\lambda'$, we have
$\der\lambda'=i\mu'\dz\bar{\mu}'+(c'\mu'+\bar{c}'\bar{\mu}')\dz\lambda'$,
with $c'=h\big(c-t_0-\partial\log(h\bar{h})\big)$. Thus, if
our CR structure admits $\mu'$ of (\ref{323a}) then
$(\lambda',\mu',\bar{\mu}')$ satisfy, qualitatively, the same structural
equations (\ref{onecr})-(\ref{psec}) as
$(\lambda,\mu,\bar{\mu})$. Therefore, in such a situation, when lifting the
CR manifold $M$ to a spacetime satisfying Einstein equations (a), we
can use $(\lambda,\mu,\bar{\mu})$ and $(\lambda',\mu',\bar{\mu}')$ on 
the same footing. We know that if we start with
$(\lambda,\mu,\bar{\mu})$, we get our conclusions
(\ref{ee1})-(\ref{ee5}). 
Similarly, using $(\lambda',\mu',\bar{\mu}')$
we get the same conclusions, with the mere change that all the variables
in (\ref{ee1})-(\ref{ee5}) have now \emph{primes}. It is easy
to get the relations between the `primed' and the `nonprimed'
variables. For us the most important is the relation between $t$
calculated for $(\lambda,\mu,\bar{\mu})$ and $t'$ calculated for
$(\lambda',\mu',\bar{\mu}')$. This, when calculated, is 
\be
t'=h(t-t_0). \label{ttt}
\ee

The hypothetic situation, in which we have $\mu'$ as in 
(\ref{323a}), is realized in practice if we have a CR
structure as in Section \ref{ree} which satisfies equations (a). For such
a structure, chosing $\mu=\der\zeta$, we get $t$ satisfying
(\ref{ee5}). Then, given such a $t$, we define
$\varphi=\mu+i\bar{t}\lambda$ as in (\ref{0psy}). Since, as we know, this
$\varphi$ satisfies $\der\varphi\dz\varphi=0$,
$\varphi\dz\bar{\varphi}\neq 0$, then by Lemma \ref{le}, we are
guaranteed an existence of $h\neq 0$ such that
$\mu'=h^{-1}(\mu+i\bar{t}\lambda)$ is closed, $\der\mu'=0$. Thus,  
passing from $(\lambda,\mu,\bar{\mu})$ to
$(\lambda'=|h|^{-2}\lambda,\mu',\bar{\mu}')$, as in (\ref{lap}), we
must use $t_0=t$ (compare the present $\mu'$ with this of (\ref{323a})). This means that, after transforming all the variables 
appearing in (\ref{ee1})-(\ref{ee5}) to their primed counterparts, we
get, in particular,  
$$t'=h(t-t_0)=h(t-t)= 0,\quad\quad\emph{everywhere}.$$  
This shows that even if $t$ is not identically
zero, including the case when it is zero at a point and nonzero off
this point, we can transform $t$ to \emph{zero everywhere} by an appropriate
choice of the adapted coframe. 
\begin{remark}
That $t$ may be gauged to zero everywhere was
\emph{subconciously} known to physicists, and used by them \cite{KSMH,RT}, in
their derivations of the maximally reduced system of equations for the
\emph{algebraically special Einstein metrics}. Actually they have never
encountered our variable $t$, since at the very beginning of their
considerations, they used a very specific choice of the adapted
coframe, that forbidded $t$ to ever appear. Being aware of this
`physicists trick' \cite{RT}, we were not gauging $t$
to zero here up to now, since nonzero $t$, if it
exists, provides us with a second CR function. However, if $t$ does not
give us a second CR function at the point around which we want to
embedd our CR manifold (because, for example, it is \emph{vanishing}
at this point), we use the argument above to gauge $t$ to \emph{zero everywhere}. In
this way we proved that the two cases: $t\equiv 0$, and $t=0$ at a
point, differ by the choice of an adapted coframe. 
\end{remark}

We summarize in the following corollary.  
\begin{corollary}
If a strictly pseudoconvex CR structure can be lifted to a spacetime
which satisfies Einstein's equations (a) then, without loss of
generality, we may assume that the variable $t$ is identically equal to zero in all of the reduced Einstein 
equations (\ref{ee1}),
(\ref{ee2}), (\ref{ee3}), (\ref{ee4}), (\ref{ee5}), (\ref{ee6}),
(\ref{ee7}), (\ref{ee8}) and in the equation $R_{33}=0$.  
\end{corollary}

In accordance with this corollary, we now put $t=0$ everywhere and
look for the second CR function in terms of other quantities than
$t$. In the rest of the paper we will frequently use the following
crucial lemma.
\begin{lemma}\label{lok}
Suppose that a strictly pseudoconvex CR manifold $M$ is represented by
forms $\lambda$ and $\mu=\der\zeta$, with $\mu\dz\bar{\mu}\neq 0$, 
as in (\ref{onecr})-(\ref{psec}). Then, if in addition, $M$ admits a
solution to the equation
\be
\partial_0\partial\bar{\eta}'=0,\label{lee}
\ee
for a complex valued function $\eta'$ on $M$ such that
\be
\partial_0\bar{\eta}'\neq 0,\label{lee1}
\ee 
then $(M,(\lambda,\mu))$ is locally embeddable. Here, as always, 
the operators $(\partial,\bar{\partial},\partial_0)$ are dual to $(\mu,\bar{\mu},\lambda)$.
\end{lemma}    
{\bf Proof} Since $\partial_0$ is a \emph{real} operator, we can always
find a real function $s$ on $M$ such that locally
$(\zeta,\bar{\zeta},s)$ are coordinates on $M$ and
$\partial_0=\frac{\partial}{\partial s}$. Given a solution $\eta'$ to
(\ref{lee}) we calculate $\bar{z}=-\partial\bar{\eta}'$ obtaining
$z=z(\zeta,\bar{\zeta},s)$. We restrict this function to the
hypersurface $s=0$ in $M$ getting $z_0=z(\zeta,\bar{\zeta},0)$. 
We search for a function $\omega_0=\omega_0(\zeta,\bar{\zeta})$ on 
$s=0$ in $M$, such that
$$\partial\bar{\omega}_0-\bar{z}_0=0.$$ This equation, as
the conjugate of the inhomogeneous CR equation in the complex plane, can
always be locally solved for $\omega_0$. Given such an $\omega_0$ on $s=0$ we
extend it to a complex valued function $\omega$ in $M$ by the
requirement that it is constant
along $s$, $$\frac{\partial}{\partial s}\omega\equiv
0,\quad\quad \omega_{|s=0}\equiv\omega_0.$$ Now we define $$\eta=\eta'+\omega.$$ We have 
$\partial_0\partial\bar{\eta}=\partial_0\partial\bar{\eta}'+\partial_0\partial\bar{\omega},$
and since $\partial_0\partial\bar{\eta}'\equiv 0$ and the commutator
\be
\partial_0\partial-\partial\partial_0=c\partial_0,\label{commu}\ee we get
$\partial_0\partial\bar{\eta}=\partial\partial_0\bar{\omega}+c\partial_0\bar{\omega}\equiv
0.$
Thus our complex function $\eta$ satisfies
\be
\frac{\partial}{\partial s}(\partial\bar{\eta})=0
\label{sds}\ee everywhere. Since
we have chosen $\omega$ so that 
$$(\partial\bar{\eta})_{|s=0}=(\partial\bar{\eta}'+\partial\bar{\omega})_{|s=0}=-\bar{z}_0+\partial\bar{\omega}_0=0,$$ 
then equation (\ref{sds}), considered as a differential equation 
for the unknown $\partial\bar{\eta}$,  
satisfies the initial condition $(\partial\bar{\eta})_{|s=0}=0$. Thus, 
$\partial\bar{\eta}$ must vanish everywhere,
$\partial\bar{\eta}\equiv 0$. This proves that if 
$\eta'$ satisfies (\ref{lee}) we have a \emph{CR function} $\eta$
associated with it. Moreover, because of our assumption (\ref{lee1})
we have  
$$\der\eta\dz\der\zeta\dz\bar{\mu}=\der\eta\dz\mu\dz\bar{\mu}=(\der\eta'+\der\omega)\dz\mu\dz\bar{\mu}=(\partial_0\eta')\lambda\dz\mu\dz\bar{\mu}\neq
 0.$$ This in particular means that $\der\zeta\dz\der\eta\neq 0$. Thus 
the CR functions $\zeta$ and $\eta$ are
 independent and as such they provide a local embedding of the CR
 manifold $M$ in $\bbC^2$. This finishes the proof of the Lemma.

The rest of the paper uses this Lemma, under various further
assumptions about the lifted spacetime, to produce a new
CR function which, together with the $\zeta$ of (\ref{koty}), provides the
embedding of the CR manifold. 
\subsubsection{Existence of a null Maxwell field aligned with the
  congruence}\label{maxz}
As a warm up we start with a CR manifold $M$ and assume 
it lifts to a spacetime that merely satisfies Einstein equations
(a). As we know, in such a case, we automatically have the CR function
$\zeta$ of (\ref{koty}) which can be used to choose the forms
$\lambda$ and $\mu$ as in (\ref{onecr})-(\ref{psec}). Next we add the 
assumption about the corresponding spacetime metric (\ref{mde}). We
will assume for a while that the lifted spacetime $({\mathcal M},g)$
admits a \emph{null Maxwell field} which is \emph{aligned} with the
congruence of null geodesics corresponding to $(M,(\lambda,\mu))$. The
terms in itallics mean the following:
\begin{itemize}
\item In any oriented spacetime $({\mathcal M},g)$ a \emph{Maxwell field} is
  a real 2-form $F$ such that $\der F=\der*F=0$, where $*$ is the Hodge
  star operator associated with the metric $g$.
\item Every real 2-form $F$ in spacetime defines a complex 
  2-form ${\mathcal F}=F+i*F$. This is antiselfdual, i.e. by
  definition,  
it satisfies $*{\mathcal
  F}=-i{\mathcal F}$. Also the converse is true. Every complex
antiselfdual 2-form
  $\mathcal F$ defines a real
  2-form $F$, via $F={\rm Re}{\mathcal F}$. The so defined $F$ has the
  property that ${\rm
  Im}{\mathcal F}=*F$. 
\item Thus Maxwell fields are in one to one correspondence with
  closed antiselfdual 
complex 2-forms $\mathcal F$ in ${\mathcal M}$. From now on we will identify Maxwell fields
  with such $\mathcal F$s.
\item A nonzero Maxwell field $\mathcal F$ is called \emph{null} iff
  ${\mathcal F}\dz{\mathcal F}\equiv 0$. Thus a null Maxwell is
  algebraically special.
\item An example of a null Maxwell field is given by a plane
  electromagnetic wave, in which the electric field ${\bf E}$ and the magnetic 
field ${\bf B}$ are orthogonal to each other ${\bf EB}=0$ and have
  equal length ${\bf E}^2-{\bf B}^2=0$. In this case $F=\der t\dz({\bf
  E}\der{\bf r})+\tfrac12\der{\bf r}\dz({\bf B}\times\der{\bf r})$ and
  $*F$ is obtained from $F$ by the replacement ${\bf E}\to{\bf B}$ and
  ${\bf B}\to-{\bf E}$.
\item A nontrivial example is due to Ivor Robinson \cite{aat}. We
  present it here, due to its influence on the entire
  subject: 

In Minkowski spacetime ${\mathcal M}=\bbR^4$, with the
  metric $g=2(\der u'\der r+\der \zeta\der\bar{\zeta})$, consider the
  following (complex) change of variables:
$$u'=u-r z\bar{z},\quad
  \zeta=(r+i)z,\quad\bar{\zeta}=(r-i)\bar{z}.$$ After this
  transformation the metric is $g=2\big(\lambda\der r+2(r^2+1)\mu\bar{\mu}\big)$, 
with $\lambda=\der u+i(\bar{z}\der z-z\der\bar{z})$ and $\mu=\der z$. Now
consider an antiselfdual 2-form ${\mathcal
  F}=f\lambda\dz\mu$ with a nonvanishing (sufficiently smooth) complex valued function $f$ in $\mathcal M$. It is obviously
\emph{null}, and it defines a \emph{Maxwell} field, i.e. it satisfies $\der{\mathcal F}=0$, if and only
if 
\begin{itemize}
\item $f=f(u,r,z,\bar{z})$ is independent of the real coordinate $r$,  
$f_r=0$, and 
\item $f$ satisfies the \emph{linear} PDE:  
$\frac{\partial f}{\partial\bar{z}}+i z\frac{\partial f}{\partial
  u}=0.$
\end{itemize}
The beauty of this example is that  
$$\frac{\partial}{\partial\bar{z}}+i z\frac{\partial}{\partial
  u}$$ is the Hans Lewy operator \cite{Le2,aat}. 
\item Null Maxwell fields are radiative in a similar sense as the
  algebraically special gravitational fields. Far from the sources the
  leading term of the field strength $F$ behaves as 
$F=\frac{\rm Null}{r}+O(\frac{1}{r^2})$, as $r\to\infty$.
\item A null Maxwell field is always of the form ${\mathcal
  F}=\tilde{f}\tilde{\lambda}\dz\tilde{\mu}$, with some real 1-form 
$\tilde{\lambda}$, some complex 1-form $\tilde{\mu}$ and a complex 
function $\tilde{f}$ on $\mathcal M$.
\item One of the implications of the Robinson theorem \cite{robin} is
  that if the spacetime ${\mathcal M}$ admits a 
null Maxwell field ${\mathcal
  F}=\tilde{f}\tilde{\lambda}\dz\tilde{\mu}$, then it is locally a
  product ${\mathcal M}=\bbR\times M$, with $M$ being a CR
  manifold. The CR structure on $M$ is induced by the forms $\lambda$,
  $\mu$ on $M$ such that $\tilde{\lambda}=\pi^*(\lambda)$ and
  $\tilde{\mu}=\pi^*(\mu)$, with $\pi:{\mathcal M}\to M$ being the
  projection which forgets about the $\bbR$ factor in $\mathcal M$. 
\item Since a null Maxwell field in the spacetime 
induces the CR structure as above, then the congruence in $\mathcal M$,
being tangent to the $\bbR$ factor, is null geodesic and shearfree. 
\item Now the construction in the other direction can be
  attempted. Given a 3-dimensional CR structure $(M,(\lambda,\mu))$
  one considers its lift to the spacetime ${\mathcal M}=\bbR\times M$, 
which is then naturally equipped with a null congruence of shearfree
  geodesics tangent to the $\bbR$ factor. Then the null Maxwell field
  ${\mathcal F}=f\lambda\dz\mu$ is called \emph{aligned} with this congruence.  \end{itemize}     
Thus let us assume that in addition to the strict pseudoconvexity,
  and to the assumption that the lifted spacetime satisfies equations (a), we 
  also have a nonvanishing null Maxwell field aligned with the
  congruence associated to CR structure $(M,(\lambda,\mu))$. This
  additional assumption will
  play the role of our equation ($2^{\rm nd}CR)$. Assuming this, 
we are guaranteed the existence of a
  complex function $f$ on $M$ such that $$\der(f\lambda\dz\mu)=0.$$
  For $(M,(\lambda,\mu))$ as in (\ref{onecr})-(\ref{psec}) we easily
  check 
  that this equation is equivalent to
  $$(\bar{\partial}f+\bar{c}f)\lambda\dz\mu\dz\bar{\mu}=0.$$ This in
  turn is equivalent to a single complex equation \cite{Taf}
\be
\partial \bar{f}+c\bar{f}=0\label{nbm}\ee
on $M$, which is our new ($2^{\rm nd}CR$). If we now introduce a
function $\bar{\eta}'$ related to $f$ by 
\be
\partial_0\eta'=f,\label{feta}\ee
we have $\partial_0\bar{\eta}'\neq 0$, since otherwise
the Maxwell field would vanish. Moreover, inserting the definition
  (\ref{feta}) in the Maxwell equation (\ref{nbm}), we see that if
  a nonvanishing $f$ satisfies (\ref{nbm}), then the corresponding $\eta'$ satisfies
\be
\partial\partial_0\bar{\eta}'+c\partial_0\bar{\eta}'=0.\label{cruc}
\ee
Using the commutator (\ref{commu}) we conclude that this equation is
finally equivalent to $\partial_0\partial\bar{\eta}'=0$. Since, as we
have already noticed $\partial_0\bar{\eta}'\neq 0$, our present
  $\eta'$ satisfies all the assumptions of Lemma \ref{lok}. Using it
  we define a CR function $\eta$ which is independent of $\zeta$. Thus we
  have proven Theorem \ref{maxc} in the direction
  ${\rm (i)\&(ii)}\Rightarrow{\rm embeddability}$. 

To get the converse
  we do as follows. Assuming embeddability, we have two independent CR
  functions. Let us choose one, say $\zeta$, such that
  $\der\zeta\dz\der\bar{\zeta}\neq 0$. Then, using $\zeta$, 
we construct a spacetime whose Ricci tensor satisfies equations
  (a), as in the proof of Theorem \ref{fistcr}. After achieving this we, 
in particular, have ${\rm embeddability}\Rightarrow {\rm (i)}$. In
  addition we have $\mu=\der\zeta$ and $\lambda$ of (\ref{psec}). 
Now we take an independent CR function, say $\eta$. Because of the
  independence condition $\der\zeta\dz\der\eta\neq 0$, we have that
  $\partial_0\eta\neq 0$. Then we define $f=\partial_0\eta\neq 0$ and
  observe that ${\mathcal F}=f\lambda\dz\mu$ satsifies $\der {\mathcal
  F}=0$. This provides us with a nontrivial null aligned Maxwell field,
  proving that ${\rm embeddability}\Rightarrow {\rm (ii)}$. Thus
  Theorem \ref{maxc} is proven.

This completes our discussion of the existence of a null Maxwell
field in the spacetime. We mention however that Trautman
  \cite{Tr} has conjectured that Theorem \ref{maxc} remains valid
  without condition (i).

\subsubsection{Petrov type II or D} We now return to the pure Einstein
situation, in which we have a strictly pseudoconvex CR manifold $M$
whose lifted spacetime $({\mathcal M},g)$ satisfies Einstein equations
(a). We work in the gauge $t\equiv 0$ and we impose
further Einstein equations on the lifted spacetime. From now on we
will always assume that the lifted spacetime satsifies Einstein's equation
(a), (b) and one of the equations (c), namely $R_{13}=0$. These,
according to Remark \ref{rty} are equivalent to $Ric(g)=\Lambda
g+\Phi\lambda\otimes\lambda$. As a consequence we
are guaranteed the existence of a complex function $m$ on the CR
manifold $M$ such that 
\be\partial m+3cm=0\label{eee8}\ee (compare with (\ref{ee8})
assuming $t\equiv 0$). 

Einstein's equations (a), (b) and $R_{13}=0$ do not guarantee that 
$m$ is nonvanishing. This shall be
\emph{assumed}, and the equation (\ref{eee8}) with $m\neq 0$ will
be our new ($2^{\rm nd}CR$). 

The assumption about the existence of a
\emph{nonvanishing} $m$ has a clear spacetime meaning. This is due to 
equation (\ref{tyd}). It says that $m\neq 0$ at a point 
if and only if the spacetime metric at this point is of no
more special Petrov type than $II$ or $D$. So let us assume that a
strictly pseudoconvex 3-dimensional CR manifold $M$ lifts 
to a spacetime of Petrov type $II$ or $D$, but no more algebraically
special, which in addition satisfies Einstein's equations (a), (b) and
$R_{13}=0$. Having assumed this we replace equation (\ref{eee8}) with
an equation for the complex function $\eta'$ such that \cite{lnt} 
\be
m=(\partial_0\bar{\eta}')^3.\label{eee9}\ee
By our assumption about the Petrov type we have
$$\partial_0\bar{\eta}'\neq 0.$$ Moreover, inserting (\ref{eee9}) 
into (\ref{eee8}), after the trivial simplification which
uses the assumed $\partial_0\bar{\eta}'\neq 0$, we
get $\partial\partial_0\bar{\eta}'+c\partial_0\bar{\eta}'=0.$ This is
again equation (\ref{cruc}) for $\eta'$ satisfying
$\partial_0\bar{\eta}'\neq 0$. This means that the argument from the
previous subsection applies, and using Lemma \ref{lok}, we can modify 
$\eta'$ to $\eta$ being the second CR function on $M$. This proves
the following theorem. 
\begin{theorem}\label{sxc}
Assume that a strictly pseudoconvex CR manifold $M$ (i) admits a lift to a spacetime
satisfying Einstein's equations $Ric(g)=\Lambda
g+\Phi\lambda\otimes\lambda$, and (ii) has Petrov type $II$ or
$D$, but no more special. 
Then such a CR structure is locally CR embeddable.
\end{theorem}
At this stage it is worthwhile to note that if we have $m$ satisfying
(\ref{ee8}) we can use its associated $\eta'$ to define $f$ by formula
(\ref{feta}). Then we can use this $f$ to define a Maxwell 2-form 
${\mathcal F}$ by ${\mathcal F}=f\lambda\dz\mu$. Due to 
(\ref{cruc}) this $\mathcal F$ satisfies
the Maxwell equations $\der{\mathcal F}=0$. Thus the lifted spacetime
of Theorem \ref{sxc} admits an aligned null Maxwell field.   

Also the converse to Theorem \ref{sxc}, namely Theorem
\ref{petovt}, can now be proven, using a similar argument. Indeed,
given an embeddable strictly pseudoconvex CR manifold $M$, we choose
one of its CR functions $\zeta$ such that $\der\zeta\dz\der\bar{\zeta}\neq
0$ to define $\mu=\der\zeta$ and $\lambda$ satisfying (\ref{psec}). 
Then we take $t\equiv 0$ and $s\equiv 0$. To construct an Einstein
spacetime satisfying $Ric(g)=\Lambda
g+\Phi\lambda\otimes\lambda$ we need first to find a complex 
function $m$ such that equation (\ref{ee8}) with $t\equiv 0$ 
is satisfied. We can do it
in two ways. Either we choose $m\equiv 0$, or we can prove that we can
find $m\neq 0$ satisfying (\ref{ee8}).

Let us first consider the second possibility. Since our real analytic 
CR manifold is
locally embeddable \cite{ah} we are guaranteed that a second CR function $\eta$,
independent of $\zeta$ exists on $M$. Thus we have $\partial_0\eta\neq
0$. Then we define $m=(\partial_0\eta)^3$, which obviously does not
vanish. Because $\eta$ satisfies the 
tangential CR equation, we easily get that our $m$
satsifies $\partial m+3cm=0$. After determining $\eta$ we must solve
the last of the reduced Einstein equations $Ric(g)=\Lambda
g+\Phi\lambda\otimes\lambda$, namely equation (\ref{ee7}) for
$p$. This is a \emph{real} equation for a \emph{real} function on the 
CR manifold $M$. Looking at (\ref{ee7}) we see that if we specify a
real constant $\Lambda$, the only unknown
in this equation is $p$, since now $c$, $t\equiv 0$ and $m$, as 
well as $\partial$ and $\bar{\partial}$ are just specified. 
For every fixed constant $\Lambda$, this equation is a 
second order real PDE in 3-dimensions, with a quite well behaved 
nonlinear part. We do not know for sure about its solvability unless we 
assume real analyticity in the variables appearing in it. To have this it
is enough to assume that the CR manifold $M$ is \emph{real
  analytic}. Then we can always find a local solution for $p$. Inserting this $p$, and $m$, $t\equiv 0$, $s\equiv
0$, $c$ into (\ref{mde}), with functions $P$, $W$ and $H$ as in
(\ref{ee1}), (\ref{ee2}), (\ref{ee3}), (\ref{ee4}), (\ref{ee6}) we
define a metric $g$. This satisfies $Ric(g)=\Lambda
g+\Phi\lambda\otimes\lambda$, with our fixed $\Lambda$ and some real 
function $\Phi$, which is determined by all our choices. Due to the
comitted choice of
$\eta$, we have $m\neq 0$, therefore the lifted spacetime is of
Petrov type $II$ or $D$. This proves Theorem
\ref{petovt}. 

\begin{remark}\label{rekm}
We strongly believe that equation (\ref{ee7}) with arbitrary
sufficiently smooth functions $t$, $c$, $m$ and arbitrary real
constant $\Lambda$, has a local nonvanishing solution for $p$ on any
sufficiently smooth CR manifold $M$, and that it could be proved by
standard methods. If this is true then we could
replace Theorem \ref{petovt} with a stronger one, in which the term `real
analytic' would be replaced by `sufficiently smooth embeddable'.  
\end{remark}

The authors are unaware of a precise reference to the literature in
which the existence of nonzero solutions to (\ref{ee7}), without the 
real analyticity assumption, is
proved. In the rest of this section we will \emph{assume} that it
is true.\\

Let us now discuss the first possibility mentioned
above. Actually, instead of using the second CR function, we could
have chosen $m\equiv 0$ in addition to $t\equiv 0$ and $s\equiv
0$. Then inserting these functions into the equation (\ref{ee7}) for
$p$ and fixing a constant $\Lambda$, we conclude that it admits a local solution. Thus defining the
metric as before we again get a spacetime with $Ric(g)=\Lambda
g+\Phi\lambda\otimes\lambda$. There is however an important difference
between this situation and the one considered before. The spacetime
now has $m\equiv 0$, so that it is of Petrov $III$ or its
specialization. Even more important is the fact that in constructing
the metric now we \emph{did not use the second CR function}. Thus,
modulo our current assumption, we have the following corollary. 
\begin{corollary}
Every sufficiently smooth 
strictly pseudoconvex 3-dimensional CR structure which admits one
CR function $\zeta$ such that $\der\zeta\dz\der\bar{\zeta}\neq 0$
has a lift to a spacetime which satisfies Einstein equations $Ric(g)=\Lambda
g+\Phi\lambda\otimes\lambda$ and which is of Petrov type $III$, or its
specializations $N$ or $0$. 
\end{corollary}       
This means that for the price of generality in the Petrov type, we may 
replace the embeddability assumption in Theorem \ref{petovt}, by a \emph{weaker} assumption about the mere
existence of one CR function, and still get the Einstein condition $Ric(g)=\Lambda
g+\Phi\lambda\otimes\lambda$ for the lift. 
\subsubsection{$Ric(g)=0$ and Petrov type III }
If we assume that our strictly pseudoconvex CR structure $(M,(\lambda,\mu))$ has a lift
to a spacetime $({\mathcal M},g)$ satisfying Einstein's
equations $Ric(g)=\Lambda g+\Phi\lambda\otimes\lambda$, which in addition, is of Petrov type 
$III$ or its specializations, then \emph{without further assumptions} 
about $({\mathcal M},g)$ we are \emph{unable} to
produce the second CR function for $M$. Of course, to get the embeddability of
$M$, we may assume that our spacetime admits an aligned null Maxwell
field and then use Theorem \ref{maxc}. But, if we lack a Maxwell field detector, we need
to invent a new ($2^{\rm nd}CR$) equation that guarantees the
existence of a second CR function $\eta$. This can be done by imposing
more special restrictions on $Ric(g)$, as we will do in this section.

So now we assume that the lifted spacetime of our CR manifold
$(M,(\lambda,\mu))$ is of Petrov type $III$ or more special, and that
it satsifies Einstein's equations (a), (b) and the first of
equations (c), namely $R_{13}=0$. We work in the gauge $$t\equiv 0$$ and, 
due to our assumption about the Petrov type, we have $$m\equiv 0.$$ Then, 
guided by the theory of exact solutions of Einstein equations we 
introduce a function 
\be
I=\partial(\partial\log p+c)+(\partial\log p+c)^2.
\ee 
This enable us to significantly simplify the formulae for the last
componenent of the Ricci tensor and the Weyl scalar coefficient
$\Psi_3$. These are given in the following proposition.  
\begin{proposition}
If the lifted spacetime satsifies the Einstein equations (\ref{ee1}),
(\ref{ee2}), (\ref{ee3}), (\ref{ee4}), (\ref{ee5}), (\ref{ee6}),
(\ref{ee7}), (\ref{ee8}) then the Ricci tensor component $R_{33}$ is
given by
\be
R_{33}=8\frac{\cos^4(\frac{r+s}{2})}{p^4}~\Big(\partial +2c
\Big)(p^2\partial \bar{I})+O(\Lambda)\label{r33}
\ee
and the Weyl scalar $\Psi_3$ is given by
\be
\Psi_3=2i~\partial
  \bar{I}~\frac{{\rm e}^{\frac{i(r+s)}{2}}}{p^2}\cos^3(\frac{r+s}{2})+O(\Lambda),\label{psi3}
\ee
as $\Lambda\to0$.
\end{proposition}
\begin{remark}
The omitted $O(\Lambda)$ term in $R_{33}$ reads 
\begin{eqnarray*}
&&O(\Lambda)=-8\Lambda\cos^4(\frac{r+s}{2})~\times\\
&&\Big(\frac43\Lambda
p^2+6(\bar{c}\partial+c\bar{\partial})\log p+12\partial\log
p\bar{\partial}\log p+3c\bar{c}-\tfrac12(\partial\bar{c}+\bar{\partial}c)-2i\partial_0\log p\Big),\end{eqnarray*}
and the omitted $\Lambda$ term in $\Psi_3$ is
$$O(\Lambda)=-4i~\Lambda~(2\bar{\partial}\log p+\bar{c})~{\rm e}^{\frac{i(r+s)}{2}}\cos^3(\frac{r+s}{2}).$$
We note that the $O(\Lambda)$ term in $R_{33}$ is
\emph{complex}. It includes the purely imaginary $-2i\partial_0\log p$. Thus the
first term in $R_{33}$ is \emph{also complex}, since $R_{33}$ is 
\emph{real}. This means that the first term in $R_{33}$ includes a
purely imaginary $\Lambda$ term which cancels $-2i\partial_0\log
p$. If $\Lambda=0$ the first term in $R_{33}$ becomes real, and
$R_{33}$ becomes real as it should be. 
\end{remark}
The appearance of the unwanted 
$O(\Lambda)$ terms in (\ref{r33}) and (\ref{psi3}) forces us to
assume that $\Lambda=0$. So in our search for the second CR function
we will assume from now on that the lifted spacetime has vanishing
cosmological constant
$$\Lambda=0.$$
Then, if we in addidtion assume that the lifted spacetime is
Ricci flat, we may easily use the function $\partial\bar{I}$ to construct
the second CR function. 

Let us thus assume that the lifted spacetime has $\Lambda=0$ and
$R_{33}=0$ everywhere, and that in addition it is of strictly Petrov
type $III$. This last assumption means that $\partial \bar{I}\neq
0$. Moreover, since $R_{33}=0$ and $\Lambda=0$ guarantees that
\be
(\partial+2c)(p^2\partial \bar{I})=0,\label{r33a}\ee
we may use our standard trick of considering $\eta'$ related to $I$
via:
$$p^2\partial \bar{I}=(\partial_0\bar{\eta}')^2.$$  
Inserting this into (\ref{r33a}) and utilising the assumption
$\partial_0\bar{\eta}'\neq 0$ about the Petrov type, we again obtain 
$\partial_0\partial\bar{\eta}'=0$, which is enough to conclude that the
following theorem is true:
\begin{theorem}\label{tt3}
Assume that a strictly pseudoconvex CR manifold $M$ admits a lift to
the spacetime satisfying (i) Einstein's equations $Ric(g)=0$, and
(ii) having 
Petrov type $III$, but no more special. 
Then such a CR structure is locally embeddable.
\end{theorem} 
Of course, as in the end of the last subsection we can now use our
second CR function, to construct an aligned null Maxwell field in our
Ricci flat spacetime of type $III$.
\subsubsection{Petrov type N}
Staying in the gauge $t\equiv 0$, with the cosmological constant set
to $\Lambda =0$, an assumption that our spacetime is
of type $N$ means that $m=0$ and $\partial \bar{I}=0$ everywhere 
(see (\ref{tyd}) and (\ref{psi3})). In the context of our search for
the second CR function this is a very fortunate Petrov type. Indeed,
assuming type $N$ we have $$\partial \bar{I}\equiv 0,$$
which not only \emph{implies} that $R_{33}=0$, but also implies that
$I$ is a CR function! The only question is if this CR function is
\emph{independent} of $\zeta$. For this we need 
\be\partial_0\bar{I}\neq 0.\label{indN}\ee
To conclude the independence we
calculate the last Weyl scalar $\Psi_4$. Assuming that $t\equiv 0$, 
$\Lambda=0$, $m\equiv 0$ and $\partial \bar{I}\equiv 0$ we get:
$$\Psi_4=2i~\partial_0 \bar{I}~\frac{{\rm
  e}^{-\frac{i(r+s)}{2}}}{p^2}\cos^3(\frac{r+s}{2}).$$
Thus the condition (\ref{indN}) for $I$ to be an independent CR
  function is equivalent to the condition on Petrov type not to be
  degenerate to the conformally flat type $0$. Thus we 
have the following theorem. 
\begin{theorem}\label{ttn}
Assume that a strictly pseudoconvex CR manifold $M$ admits a lift
to a spacetime satisfying (i) Einstein equations (a), (b) with
$\Lambda =0$ and
$R_{13}=0$ and, which in addition, (ii) has Petrov type $N$ and is
nowhere degenerate to Petrov type 0. Then such a CR structure is
locally embeddable. Moreover, in such case the
spacetime is Ricci flat.
\end{theorem} 
The remark about the existence of the aligned Maxwell field, as at the
end of the previous subsections, applies here also.
\subsubsection{Conformally flat case}
If we only know that 
among the lifted spacetimes of a strictly pseudoconvex CR structure there is a
Minkowski metric, we may proceed with our search for the second function in
the same spirit as we were doing in the previous
subsections. However, in such a case there is a simpler more elegant
geometric way of achieving our goal. This comes from
Penrose's \emph{twistor theory}.

Let us now forget about all the results from the entire Section
\ref{sefd} and assume that we are given a 3-dimensional 
CR structure $(M,(\lambda,\mu))$, not neccessarily strictly 
pseudoconvex (!), which has a lift to a \emph{conformally flat}
spacetime $\mathcal M$. We do not need the Einstein equations for 
the rest of the argument. It is known (see e.g. \cite{PenRin}) that 
the space of all null geodesics in a neighbourhood in $\mathcal M$ 
is a 5-dimensional CR manifold $N$, which is naturally 
locally CR embedded in $\bbC^3$. Actually $N$ may be identified with an open set in the following
\emph{real projective quadric} 
$$PN=\{\bbC\bbP^3\ni(W_1:W_2:W_3:W_4)~|~|W_1|^2+|W_2|^2-|W_3|^3-|W_4|^2=0\}$$ 
CR embedded in $\bbC\bbP^3$. The manifold $PN$ is called the
\emph{space of projective twistors}.

Since points of $N$ are
null geodesics in $\mathcal M$ then a congruence of null geodesics in
$\mathcal M$ is just a 3-dimensional manifold $M_N$ in $N$. Crucial to
our argument is the fact that if a congruence of null geodesics in
$\mathcal M$ is \emph{shearfree} then $M_N$ is a \emph{CR submanifold} 
\cite{nurtra} of the
5-dimensional \emph{embedded} CR manifold $N$. Thus having a
congruence of null and shearfree geodesics in $\mathcal M$ we first
are guaranteed that $M$ is locally CR embedded as a submanifold $M_N$ in
$\bbC^3$. But this implies that $M$ also has a local CR embedding in a
$\bbC^2$, see \cite{ah,hilt}. The argument is very simple:

Take a point $p\in M_N$, and define $\bbC^2$ to be the smallest complex
vector space which contains ${\rm T}_pM_N$. The local projection
$\pi:\bbC^3\to\bbC^2$ is holomorphic, hence its restriction
$\phi=\pi_{|M_N}$ is a CR map, whose image in $\bbC^2$ is the desired
CR embedding. This proves the following theorem. 
\begin{theorem}\label{tt0}
Every 3-dimensional CR manifold which has a lift to a conformally flat
spacetime is locally embeddable.
\end{theorem}

Using Theorems \ref{sxc}, \ref{tt3}, \ref{ttn}, and
\ref{tt0} we obtain Theorem \ref{petrovt}.

\end{document}